\documentclass[12pt]{article}

\usepackage{latexsym}

	\newcommand{\be}{\begin{equation}}
	\newcommand{\ee}{\end{equation}}
	\newcommand{\bea}{\begin{eqnarray}}
        \newcommand{\eea}{\end{eqnarray}}
	\newcommand{\bean}{\begin{eqnarray*}}
	\newcommand{\eean}{\end{eqnarray*}}

	\newtheorem{Thm}{Theorem}[section]
	\newenvironment{Theorem}{\begin{Thm} \sl}{\end{Thm}}

	\newtheorem{Prp}[Thm]{Proposition}
	\newenvironment{Proposition}{\begin{Prp} \sl}{\end{Prp}}

	\newtheorem{Lma}[Thm]{Lemma}
	\newenvironment{Lemma}{\begin{Lma} \sl}{\end{Lma}}

	\newtheorem{Cor}[Thm]{Corollary}
	\newenvironment{Corollary}{\begin{Cor} \sl}{\end{Cor}}

	\newtheorem{Defin}[Thm]{Definition}
	\newenvironment{Definition}{\begin{Defin} \rm}{\end{Defin}}

	\newtheorem{Examp}[Thm]{Example}

	\newtheorem{Examps}[Thm]{Example}

	\newenvironment{Proof}{{\noindent \sc Proof:}}{\hfill $\Box$ \vspace{5 mm}}

\newcommand {\mybib} [6]
{
\bibitem[#1]{#2} #3 {\sl #4}. #5, #6
}
\newcommand {\serbib} [7]
{
\bibitem[#1]{#2} #3 {\sl #4}. #5 (#6), #7
}

\newcommand{\BR}{\mathbf{R}}
\newcommand{\BZ}{\mathbf{Z}}
\newcommand{\BN}{\mathbf{N}}
\newcommand{\BC}{\mathbf{C}}
\newcommand{\BT}{\mathbf{T}}
\newcommand{\half}{\frac{1}{2}}
\newcommand{\hl}{{[ \! [}}
\newcommand{\hr}{{] \! ]}}
\newcommand{\dil}{\mathcal{D}}

	\title{Wavelets and Hilbert Modules}
	\author{Peter John Wood \\ 
School of Informatics and Engineering \\
Flinders University of South Australia \\ GPO Box 2100, Adelaide, SA, Australia 5001 \\
email: pwood@infoeng.flinders.edu.au, \\
       pjwood@myplace.net.au} 
\begin{document}

\maketitle

\begin{abstract}
A Hilbert $C^*$-module is a generalisation of a Hilbert space for which the inner product takes its values in a $C^*$-algebra instead of the complex numbers. We use the bracket product to construct some Hilbert $C^*$-modules over a group $C^*$-algebra which is generated by the group of translations associated with a wavelet. We shall investigate bracket products and their Fourier transform in the space of square integrable functions in Euclidean space. We will also show that some wavelets are associated with Hilbert $C^*$-modules over the space of essentially bounded functions over higher dimensional tori.
\vspace{5 mm}


\vspace{5 mm}

\footnotetext{{\textit{Math Subject Classifications.}} 42C40, 46L08, 43A25.}
\end{abstract} 

We shall examine in detail a construction employing Hilbert $C^*$-modules that relates to wavelet theory. A Hilbert $C^*$-module (also known as a Hilbert module or a $B$-rigged space) is in some ways similar to a Hilbert space except that the inner product takes its value in a $C^*$-algebra instead of the complex numbers.  The $C^*$-algebra valued inner product that we shall use is sometimes known as the bracket product. The Hilbert $C^*$-module described here has its linear space contained in 
$L^2(\BR^d)$. Much of the work in this paper was done during the author's PhD \cite{w1}.

The bracket product has been used in wavelet theory before, see for example \cite{bdr1}, \cite{f1} and \cite{bcmo1}. The connection between Hilbert $C^*$-module theory and wavelet theory that is being investigated here was described in a talk given by M. A. Rieffel in 1997 \cite{r3}. The material in \cite{r3} has more recently been elaborated on in two papers by J. A. Packer and M. A. Rieffel \cite{pr1, pr2}.

In \cite{pr1}, a Hilbert $C^*$-module is constructed with a linear space contained in $C^*(\BZ^d)$, and is used to study the properties of continuous filters. The Hilbert $C^*$-modules described in \cite{pr2} are similar to the ones described here. The paper \cite{pr2} contains some interesting work investigating wavelets in arbitrary projective modules over $C^*(\BZ^d)$. In this paper we shall be focusing on using Hilbert $C^*$-modules to study classical wavelet theory. The Hilbert modules that we shall construct can be thought of as a special case of the main example described in \cite{r1}. This example is constructed from a closed subgroup of a locally compact group and is also described in \cite[Example 1]{r1} and \cite[Appendix C]{rw1}.

A related Hilbert $C^*$-module construction which also uses bracket products has been developed by P. G. Casazza, M. Coco, and M. C. Lammers and is described in \cite{cala1} and \cite{cola1}. The Hilbert $C^*$-module described in \cite{cala1} and \cite{cola1} is over the $C^*$-algebra $L^\infty([0,1])$ and is used to study Gabor systems. In this paper we use the bracket product to construct Hilbert $C(\BT^d)$-modules and Hilbert $L^\infty(\BT^d)$-modules in the Fourier domain. The Hilbert $C(\BT^d)$-modules are also $C^*(\BZ^d)$-modules in the time domain.

Our construction is related to some work on frames for Hilbert $C^*$-modules by M. Frank and D. Larson in \cite{fl2, fl1, fl3}. Other results which relate wavelets to $C^*$-algebras and other operator algebras have been obtained in \cite{bj2, bj1, bj3, dl1, hl2, lr1, rt1}.

%
%

In Section \ref{s1} we shall introduce the main notation that shall be used and review some background material on wavelets. In Section \ref{s2} we shall provide the necessary background material on Hilbert $C^*$-modules.

The main construction is fleshed out in detail in Sections \ref{s3}, \ref{s4} and \ref{s5}. In Section \ref{s3} we shall consider an embedding $\theta : \BZ^d \to \BR^d$, and use it to construct a Hilbert $C^*(\BZ^d)$-module $X_\theta$, which we shall show is contained in $L^2(\BR^d)$. We shall also construct a Hilbert $C(\BT^d)$-module  $\hat{X_\theta}$, which is the image of $X_\theta$ under the Fourier transform.
In Lemma \ref{lma:rdfouhilb} we analyse the role of the Fourier transform in detail. In Lemma \ref{lma:hmnormx} we show that the Hilbert modules that we construct can be embedded in $L^2(\BR^d)$. The main result in this section is Theorem \ref{rdhmexists} which verifies that $X_\theta$ and $\hat{X_\theta}$ are full Hilbert modules, and that the Fourier transform defines a Hilbert module isomorphism.

In Section \ref{s4} we shall then incorporate the dilation by constructing a chain of Hilbert $C^*(\BZ^d)$-modules $(X_n)_{n \in \BZ}$. For each integer $n$, the Hilbert $C^*(\BZ^d)$-module $X_n$ is associated with the action of $\Gamma^n$ on $L^2(\BR^d)$. We shall construct Hilbert $C(\BT^d)$-modules $(\hat{X_n})_{n \in \BZ}$ which are the images of $(X_n)_{n \in \BZ}$ under the Fourier transform on $\BR^d$. The main result of this section in Theorem \ref{rdwavelethmexists}, which uses Theorem \ref{rdhmexists} to verify that $X_n$ and $\hat{X_n}$ are full Hilbert modules. We then prove Corollary \ref{tdadjointable}, which applies some basic Hilbert module theory to our construction. Corollary \ref{tdadjointable} demonstrates that the translations and dilation can be thought of as adjointable operators. Proposition \ref{xnthesame} is based on \cite[Proposition 1.11]{pr2}, and shows that each of the Hilbert $C^*(\BZ^d)$-modules $(X_n)_{n \in \BZ}$ share the same linear space.

In Section \ref{s5} we show how filters can be defined using the bracket product. We use this technique to show that the Shannon wavelet is not contained in the linear space corresponding to the Hilbert $C^*$-modules defined in the previous section. So we define a chain of Hilbert $L^\infty(\BT^d)$-modules $(Y_n)_{n \in \BZ}$ which are based on a construction described in \cite{cala1} and \cite{cola1}. Each $Y_n$ contains $\hat{X_n}$ as a subset (here we only work in the Fourier domain). These Hilbert modules have the property that every wavelet is contained in its linear space. The main result in this section is Theorem \ref{thm:linfty} which is used to verify that each $Y_n$ is a Hilbert $L^\infty(\BT^d)$-module.

\section{Wavelet definitions} \label{s1}

The classical definition of a wavelet in $L^2(\BR)$ is a function $\psi$ such that the family
$$
U_\psi := \left \{ 2^{-j/2} \psi(2^j x - k) \right \}_{j,k \in \BZ}
$$
is an orthonormal basis for $L^2(\BR)$. Eack element of $U_\psi$ is obtained from $\psi$ by an integer translation $\psi(x) \mapsto \psi(x-k)$	followed by a dilation $\psi(x) \mapsto 2^{j/2} \psi(2^j x)$ by a power of 2. The space of square integrable functions $L^2(\BR)$ is an example of a Hilbert space and the translations and dilations are unitary operators on this Hilbert space. This allows us to generalise this definition. 

\begin{Definition} \label{def:mrs}
A {\em multiresolution structure} on a separable infinite dimensional Hilbert space ${\mathcal H}$ is a discrete group $\Gamma$ of unitary operators on ${\mathcal H}$ and a unitary operator ${\mathcal D}$ on ${\mathcal H}$ such that ${\mathcal D}^{-1}\gamma{\mathcal D} \in \Gamma$ for all $\gamma \in \Gamma$. We call $\Gamma$ the group of {\em translations}, and ${\mathcal D}$ the {\em dilation}. 

For a multiresolution structure the set $\mathcal{D}^{-1} \Gamma \mathcal{D}$ is a subgroup of $\Gamma$. If $m$ is the index of the group $\mathcal{D}^{-1} \Gamma \mathcal{D}$ in $\Gamma$, we call $m$ the {\em index} of the multiresolution structure (recall that the index of a subgroup in a group is the number of cosets of the subgroup in the larger group). We assume that $m$ is finite.

A function $\psi \in {\mathcal H}$ is an {\em orthonormal wavelet} with a multiresolution structure $(\Gamma, {\mathcal D})$ if $\{{\mathcal D}^n(\gamma(\psi))\}_{\gamma \in \Gamma, n \in \BZ}$ is an orthonormal basis for ${\mathcal H}$. An {\em orthonormal multiwavelet} with a multiresolution structure $(\Gamma,{\mathcal D})$ is a finite set of elements $\{\psi^1, \ldots,\psi^M \}$ in ${\mathcal H}$ such that $\{{\mathcal D}^n(\gamma(\psi^i))\}_{\gamma \in \Gamma, n \in \BZ, i = 1 \ldots n}$ forms an orthonormal basis for ${\mathcal H}$.
\end{Definition}

We will abbreviate the terms ``orthonormal wavelet'' and ``orthonormal multiwavelet'' as ``wavelet'' and ``multiwavelet'', respectively.

The most important example of a multiresolution structure is given by the Hilbert space $L^2(\BR^n)$, a discrete Abelian subgroup $\Gamma$ of $\BR^n$ which induces a translation on ${\mathcal H}$ by $(\gamma f)(x) = f (x - \gamma)$ where $\gamma \in \Gamma$; and a dilation ${\mathcal D}$ given by $({\mathcal D} f)(x) = \sqrt{\det \tilde{{\mathcal D}}}f(\tilde{{\mathcal D}}x)$ where $\tilde{{\mathcal D}}$ is a mapping from $\BR^n \to \BR^n$ which also maps $\Gamma$ onto a proper subgroup of itself. Another example of a multiresolution struction is defined on the Cantor group, which is described in detail in \cite{lg1} and \cite{lg2}.

In all the cases that we will consider we will expect that ${\mathcal D}^{-1}\Gamma{\mathcal D}$ will be a proper subgroup of $\Gamma$. Define $\Delta \in \mathrm{Hom}(\Gamma, \Gamma)$ by $\Delta(\gamma) := {\mathcal D}^{-1} \gamma {\mathcal D}$. Define $\Gamma^n := {\mathcal D}^n \Gamma {\mathcal D}^{-n}$,  note that $\Gamma^{-1} = \Delta(\Gamma)$.  It is straightforward to show that $\Delta(\Gamma)$ is isomorphic to $\Gamma$ and hence for all $m,n$, $\Gamma^m \cong \Gamma^n$. We can then define a family of representations of $\Gamma$ by $\pi^n_{\gamma}(f) = \pi_{\gamma^n}(f)$ where $\gamma \in \Gamma$ and $\gamma^n \in \Gamma^n$ are mapped to each other by the isomorphism between $\Gamma$ and $\Gamma^n$.

To fix notation we define the Fourier transform of a function $f \in L^1(G)$ to be a function on the Pontryjagin dual $\hat{G}$ given by
$$
(\mathcal{F}f)(\zeta) = \hat{f}(\zeta) = \int_G f(x) \zeta(x) dx.
$$
We will make use of the {\em annihilator} $\mathrm{Ann }\Gamma$ of $\Gamma$. The annihilator of a closed subgroup $\Gamma$ of a locally compact Abelian group $G$ is the set of all $\lambda \in \hat{G}$ such that $\lambda(\gamma) = 1$ for all $\gamma \in \Gamma$. We shall sometimes use the notation $\widehat{\BR^d}$ to indicate when we are working in the Fourier domain.

\begin{Definition} \label{def:gmra}
Let $(\Gamma, {\mathcal D})$ be an multiresolution structure on a Hilbert space ${\mathcal H}$. A sequence $\{V_n\}_{n \in \BZ}$ of closed subspaces of ${\mathcal H}$ is called a {\em generalised multiresolution analysis} (GMRA) of ${\mathcal H}$ if
\begin{enumerate} 
\item $\forall n \in \BZ $, $V_n \subset V_{n+1}$
\item $ \bigcup_{n \in \BZ}V_n$ is dense in $ {\mathcal H} $ and $ \bigcap_{n
\in \BZ}V_n = \{ 0 \} $
\item $ \forall n \in \BZ, {\mathcal D}(V_n) = V_{n+1}$
\item $ V_0$ is invariant under $\Gamma$.
\end{enumerate}

A {\em multiresolution analysis} (MRA) with multiplicity $r$ is a GMRA for which there exists a set of functions $\{ \varphi^1, \ldots, \varphi^r \}$ such that $\{ \gamma(\varphi^i) \}_{\gamma \in \Gamma, i=1 \ldots r}$ is an orthonormal basis for $V_0$. We call $\{ \varphi^1 , \ldots, \varphi^r \}$ a set of {\em scaling functions}.
\end{Definition}

If we denote the orthogonal complement of $V_n$ in $V_{n+1}$ by $W_n$, it satisfies the following properties:
\begin{enumerate}
\item ${\mathcal D}(W_n) = W_{n+1}$
\item ${\mathcal H} = \bigoplus_{n \in \BZ} W_n$
\end{enumerate}
we call $W_n$ the wavelet spaces.

Let us now define what we call the standard multiresolution structure on $L^2(\BR^d)$, this is a fairly standard way of formulating wavelets on $L^2(\BR^d)$ and the main example that we shall work with. 

\begin{Definition} \label{def:standmrs}
Let $\tilde{\dil} \in M^d(\BZ)$ bs a $d \times d$ matrix with integer entries such that all of the eigenvalues of $\tilde{\dil}$ are greater than $1$. We shall call a matrix with these properties a {\em dilation matrix}. We define the {\em standard multiresolution structure} on $L^2(\BR^d)$ associated with $\tilde{\dil}$ to be the multiresolution structure with Hilbert space $L^2(\BR^d)$ and with translation group $\Gamma = \BZ^d$ (as an additive group), and dilation given by
\be \label{dildef}
(\dil f)(x) = \sqrt{m} f (\tilde{\dil} x)
\ee
for $x \in \BR^d$, $f \in L^2(\BR^d)$, and where $m \in \BN$ is the index of the multiresolution structure.
\end{Definition}

We remark that because all of the eigenvalues of $\tilde{\dil}$ are nonzero, $\tilde{\dil}$ is invertible. The following result is due to K. Gr\"ochenig and W. R. Madych (\cite[Lemma 2]{gm1}).

\begin{Lemma}
Let $\tilde{\dil}$ be as defined as in Definition \ref{def:standmrs}, and let $m$ be the index of the multiresolution structure corresponding to $\tilde{\dil}$. Then $m = |\det \tilde{\dil}|$.
\end{Lemma}

We shall sometimes work in the Fourier domain. Define
$$
\hat{\dil} = \mathcal{F} \dil \mathcal{F}^*
$$
where $\dil$ corresponds to a standard multiresolution structure. By evaluating the Fourier transform,  
one obtains that for $p \in L^2(\widehat{\BR^d})$ and $\xi \in \widehat{\BR^d}$,
\be \label{dilhatdef}
(\hat{\dil} p)(\xi) = \frac{1}{\sqrt{m}} p \left( (\tilde{\dil}^*)^{-1} (\xi) \right) .
\ee

\section{Background material on Hilbert $C^*$-modules} \label{s2}

\begin{Definition} \label{def:hmodule}
Suppose $\mathcal{A}$ is a $C^*$-algebra. An {\em inner product $\mathcal{A}$-module} is a linear space $X$ which is a right $\mathcal{A}$-module with compatible scalar multiplication: $\lambda(x a) = (\lambda x)a = x (\lambda a)$ for $x \in X$, $a \in \mathcal{A}$, $\lambda \in \BC$; and a $C^*$-algebra valued inner product $[ \; , \: ]_{\mathcal{A}} : X \times X \to \mathcal{A}$ satisfying for all $x, y, z \in X$ and $ a \in \mathcal{A}$
\begin{enumerate}
\item $[x, \alpha y + \beta z]_{\mathcal{A}} = \alpha[x,y]_{\mathcal{A}} + \beta[x,z]_{\mathcal{A}}$
\item $[x,ya]_{\mathcal{A}} = [x,y]_{\mathcal{A}} a$
\item $[x,y]^*_{\mathcal{A}} = [y,x]_{\mathcal{A}}$
\item $[x,x]_{\mathcal{A}} \geq 0$
\item $[x,x]_{\mathcal{A}} = 0 \Rightarrow x = 0$.
\end{enumerate}
\end{Definition}

It is straightforward that a $C^*$-algebra valued inner product will be conjugate-linear in the first variable and linear in the 
second variable. We can define a norm on $X$ by $\|x\| = \| [x, x]_{\mathcal{A}} \|^{{1 \over 2}}$.

\begin{Definition}
Suppose $\mathcal{A}$ is a $C^*$-algebra. An inner product $\mathcal{A}$-module which is complete with respect to its norm is called a {\em Hilbert $C^*$-module over $\mathcal{A}$}, or a {\em Hilbert $\mathcal{A}$-module}. We call a Hilbert $\mathcal{A}$-module $X$ {\em full} if $[X,X]_\mathcal{A}$ is dense in $\mathcal{A}$.
\end{Definition}

Some useful references on Hilbert modules include \cite{l1}, \cite{r1}, \cite{r2}, \cite{rw1}, and \cite{wo}.

\begin{Definition} \label{def:adjoint}
Suppose $E, F$ are Hilbert modules, with inner product in a $C^*$-algebra $A$. We define ${\mathcal L}(E,F)$ to be the set of all maps $t:E \to F$ for which there exists a map $t^*:F \to E$ such that
$$
[tx,y] = [x,t^* y]
$$
for all $x \in E, y \in F$. We call ${\mathcal L}(E,F)$ the set of {\em adjointable operators} from $E$ to $F$. We abbreviate ${\mathcal L}(E,E)$ as ${\mathcal L}(E)$. It can be shown that every element of ${\mathcal L}(E,F)$ is a bounded $A$-linear map.
\end{Definition}

\begin{Definition}
Suppose $E, F$ are Hilbert modules, with inner product in a $C^*$-algebra $A$. For $x \in E, y,z \in F$, define
$$ 
\Theta_{x,y} (z) = x \circ [y,z].
$$
It can be shown that $ \Theta_{x,y} \in {\mathcal L}(E,F) $, and $(\Theta_{x,y})^* = \Theta_{y,x}$. We define ${\mathcal K}(E,F)$ to be the closed linear subspace of ${\mathcal L}(E,F)$ spanned by $\{ \Theta_{x,y} : x \in E, y \in F \}$. We call ${\mathcal K}(E,F)$ the {\em generalised compact operators} from $E$ to $F$. ${\mathcal K}(E,F)$ is also known as the {\em imprimitivity 
algebra} of $E$ and $F$.
\end{Definition}

For $t \in \mathcal{L}(E,F)$ we define
$$
\|t\| := \sup_{\|x\| \leq 1} \| tx \| = \sup_{\|x\| \leq 1} \| [ tx,tx] \|.
$$

It can be shown that $\mathcal{L}(E)$ and $\mathcal{K}(E)$ are $C^*$-algebras with the norm defined above, and $\mathcal{K}(E)$ is an ideal in $\mathcal{L}(E)$.

\begin{Definition} \label{def:isohm}
Suppose $E, F$ are Hilbert modules, an operator $u \in \mathcal{L}(E, F)$ is called {\em unitary} if 
$$
u^* u = \mathrm{Id}, \: \: \; u u^* = \mathrm{Id} .
$$
We say that $E$ and $F$ are {\em isomorphic} if there exists a unitary $u$ contained in $\mathcal{L}(E,F)$.
\end{Definition}


The following result is very useful for constructing Hilbert $C^*$-modules from inner product pre-$C^*$-modules. It is a slightly less general version of \cite[Lemma 2.16, p15]{rw1}.

\begin{Lemma} \label{comphm}
Suppose that $\mathcal{A}_0$ is a dense $*$-subalgebra of a $C^*$-algebra $\mathcal{A}$. Suppose that $E_0$ is a right inner product $\mathcal{A}_0$-module. 
Let the linear space $E$ be the completion of $E_0$ with respect to the Hilbert module norm. Then the module action of $\mathcal{A}_0$ on $E_0$ can be extended to a module action of $\mathcal{A}$ on $E$ and the $\mathcal{A}_0$-valued inner product on $E_0$ can be extended to an $\mathcal{A}$-valued inner product on $E$ in such a way that $E$ is a right Hilbert $\mathcal{A}$-module. We call the Hilbert module $E$ the completion of the inner product module $E_0$.
\end{Lemma}
\begin{Proof}{}
See \cite[Lemma 2.16, p15]{rw1}, or \cite[Lemma 0.4.5]{w1}.
\end{Proof}

\section{Hilbert $C^*(\Gamma)$-Modules and Wavelets} \label{s3}

Consider an arbitrary embedding $\theta : \BZ^d \to \BR^d$. The proof to the following lemma is contained in the author's PhD thesis \cite[Lemma 2.1.1]{w1}.

\begin{Lemma}
\label{thetalma}
Suppose that $\theta : \BZ^d \to \BR^d$ is an embedding. There exists a unique nonsingular linear transformation $A_\theta : \BR^d \to \BR^d$ which satisfies $\theta = A_\theta \iota$, where $\iota$ is the natural embedding of $\BZ^d$ in $\BR^d$. If $\hat{\theta} : \widehat{\BR^d} \to \BT^d$ is a homomorphism which satisfies
$$
(\theta \gamma, \xi) = (\gamma, \hat{\theta} \xi)
$$
for all $\gamma \in \BZ^d$ and $\xi \in \widehat{\BR^d}$, then $\hat{\theta} = \hat{\iota} A_\theta^*$. The annihilator of $\theta(\BZ^d)$ is given by $\mathrm{Ann }\ \theta(\BZ^d) = (A^*_\theta)^{-1} \BZ^d$.
\end{Lemma}

\begin{Proof}{}
Let $\{e_i\}_{i=1}^d$ be the standard basis for $\BR^d$ and consider the set $\{ \theta (e_i )\}_{i=1}^d$. For $i,j = 1, \dots, d$ there exists $a_{ij} \in \BR$ such that $\theta(e_i) = \sum_{j=1}^d a_{ij} e_j$. If $\gamma = ( \gamma_1, \ldots, \gamma_d)$ is an arbitrary element of $\BZ^d$, then
$$
\theta(\gamma) = \sum_{i=1}^d \gamma_i \theta(e_i) = \sum_{i=1}^d \gamma_i \sum_{j=1}^d a_{ij} e_j = \sum_{j=1}^d e_j \sum_{i=1}^d \gamma_i  a_{ij} .
$$
The elements $a_{ij}$ therefore define a coordinate transformation matrix $A_\theta$. This matrix is nonsingular because otherwise there would be linear dependence among the vectors $\{ \theta (e_i) \}_{i=1}^d$, which can not occur because $\theta$ is an embedding.

Let us now show that $\hat{\theta} = \hat{\iota} A_\theta^*$. Note that because $\theta$ is an embedding it is also a group homomorphism. We calculate 
\bean
( \theta \gamma, \xi)
&=& ( A_\theta \iota \gamma, \xi)  \\
&=& ( \iota \gamma, A_\theta^* \xi) \\
&=& ( \gamma, \hat{\iota}  A_\theta^*  \xi )
\eean
so $\hat{\theta} = \hat{\iota} A_\theta^*$.

Let us now examine the annihilator of $\theta(\BZ^d)$. We calculate
\bean
\mathrm{Ann }\ \theta(\BZ^d) &=& \{ \xi \in \widehat{\BR^d} \ : \ (\theta \gamma, \xi) = 1 \mbox{ for } \gamma \in \BZ^d \} \\
&=& \{ \xi \in \BR^d \ : \ (\gamma, \hat{\iota} A_\theta^* \xi) = 1 \mbox{ for } \gamma \in \BZ^d \} \\
&=& \{ \xi \in \BR^d \ : \ A_\theta^* \xi \in \mathrm{Ann } \ \BZ^d \} \\
&=& (A_\theta^*)^{-1} \BZ^d
\eean
verifying the result.
\end{Proof}

Recall that the annihilator of $\theta(\BZ^d)$ is the set of elements $\xi$ of $\widehat{\BR^d}$ for which $\xi(x) = 1$ for any $x \in \theta(\BZ^d)$. It follows from Lemma \ref{thetalma} that $\xi \in \mathrm{Ann }\ \theta(\BZ^d)$ if and only if for all $\gamma \in \BZ^d$, $( \gamma, \hat{\theta} \xi) = 1$. This means that if $a$ is a function on $\BT^d$, then $a \circ \hat{\iota_n}$ is an $\mathrm{Ann }\ \theta(\BZ^d)$-periodic function on $\widehat{\BR^d}$, where here by $\circ$ we mean composition. Throughout this section, when we are given an embedding $\theta : \BZ^d \to \BR^d$, we shall use the notation $A_\theta$ and $\hat{\theta}$ to denote the linear transformation and dual homomorphism given by Lemma \ref{thetalma}.

\begin{Definition} \label{rdbpma}
Let $\theta$ be a embedding of $\BZ^d$ in $\BR^d$. We define a representation $\pi^\theta$ of $\BZ^d$ on $L^2(\BR^d)$ by
\be 
(\pi_\gamma^\theta f)(x) := f(x - \theta(\gamma)) \ \ \ x \in \BR^d, \ \gamma \in \BZ^d, \ f \in L^2(\BR^d).
\ee
Note that if $f \in C_c(\BR^d)$, then $\pi^\theta_\gamma f \in C_c(\BR^d)$. For $f,g \in C_c(\BR^d)$, $\gamma \in \BZ^d$, we define the {\em bracket product associated with $\theta$} to be the function on $\BZ^d$ given by
\be \label{rdbraceq}
[f,g]_\theta (\gamma) := \int_{\BR^d} \overline{f(x - \theta(\gamma))} g(x) dx .
\ee
For $f \in C_c(\BR^d)$, $a \in C_c(\BZ^d)$, we define the {\em module action associated with $\theta$} to be the function on $\BR^d$ given by
\be \label{rdbracmaeq-2}
(f \circ_\theta a) := \sum_{\gamma \in \BZ^d}a(\gamma)f(x - \theta(\gamma)) = \sum_{\gamma \in \BZ^d}a(\gamma)\pi^\theta_\gamma(f).
\ee
\end{Definition}

Note that it follows immediately from the above definition that for $f,g \in C_c(\BR^d)$,
\bea
\label{rdbpmarem1} [f,g]_\theta (\gamma) &=& \int_{\BR^d} \overline{(\pi^\theta_\gamma f)} g(x) dx \\
\label{rdbpmarem2} &=& (f^* * g)(\theta(\gamma)) 
\eea
where in (\ref{rdbpmarem1}) we embed $C_c(\BR^d)$ in $L^2(\BR^d)$.

We now summarise some of the properties of the above definitions.

\begin{Proposition} \label{bptheorem}
Let $\theta$ be an embedding of $\BZ^d$ in $\BR^d$. The bracket product $[ \: , \;]_\theta$ and associated module action $\circ_\theta$ have the following properties:
\begin{enumerate}
\item It is the case that $[f,g]_\theta \in C_c(\BZ^d)$ when $f,g \in C_c(\BR^d)$.
\item If $f \in L^2(\BR^d)$, $a \in l^1(\BZ^d)$, then the function $f \circ_\theta a$ given almost everywhere by (\ref{rdbracmaeq-2}) is measurable and contained in $L^2(\BR^d)$.
\item If $f, g \in L^2(\BR^d)$ and $a,b \in l^1(\BZ^d)$, then
\bea
\label{eqrmod1} (f + g) \circ_\theta a &=& f \circ_\theta a + g \circ_\theta a, \\
\label{eqrmod2} (f \circ_\theta a) \circ_\theta b &=& f \circ_\theta (a * b), \\
\label{eqrmod3} f \circ_\theta (a + b) &=& f \circ_\theta a + f \circ_\theta b.
\eea
It therefore follows that the module action associated with the bracket product makes $L^2(\BR^d)$ into a right $l^1(\BZ^d)$-module (where $l^1(\BZ^d)$ has the convolution product) and also makes $C_c(\BR^d)$ into a right $C_c(\BZ^d)$-module.
\item If $f,g \in L^2(\BR^d)$, then $[f,g]_\theta$ is contained in $C_0(\BZ^d)$.
\item Suppose that $f,g, h \in L^2(\BR^d)$, $a \in l^1(\BZ^d)$, $ \alpha, \beta \in \BC$, and $\theta$ is an embedding of $\BZ^d$ in $\BR^d$, then
\begin{enumerate}
\item $[f, \alpha g + \beta h]_\theta = \alpha [f,g]_\theta + \beta [f,h]_\theta$;
\item $[f, g \circ_\theta a ]_\theta = [f,g]_\theta * a$, where the binary operation $*$ is convolution on $l^1(\BZ^d)$;
\item $[f,g]_\theta^* = [g,f]_\theta$, where the unary operation $\cdot^*$ is involution on $l^1(\BZ^d)$.
\end{enumerate}
\end{enumerate}
\end{Proposition}

\begin{Proof}{}
\begin{enumerate}
\item It is true that $f \circ_\theta a \in C_c(\BR^d)$ because the set of all $w \in \BR^d$ such that there exists $\gamma \in \BZ^d$ satisying $a(\gamma)f(x - \theta(\gamma)) \neq 0$ is compact. It is the case that $[f,g]_\theta \in C_c(\BZ^d)$ because the set of $\gamma \in \BZ^d$ such that there exists $x \in \BR^d$ satisfying $\overline{f(x-\theta(\gamma))}g(x) \neq 0$ is finite.
\item Suppose that $f \in L^2(\BR^d)$ and $a \in l^1(\BZ^d)$, and let $S \subset \BZ^d$ be a finite subset. In this case the sum $\sum_{\gamma \in S} a(\gamma) f(x - \gamma)$ is a finite sum and for $x \in \BR^d$ we calculate
\bean
\left\| \sum_{\gamma \in S} a(\gamma) f(x- \theta(\gamma)) \right\|_2 
&=& \sqrt{ \int_{\BR^d} | \sum_{\gamma \in \BZ^d} a(\gamma) f(x - \theta(\gamma)) |^2 dx } \\
&\leq& | \sum_{\gamma \in \BZ^d} a(\gamma) | \sqrt{\int_{\BR^d} |f (x- \theta(\gamma)) |^2 dx } \\
&\leq& \sum_{\gamma \in \BZ^d} | a (\gamma)| \sqrt{\int_{\BR^d} |f (x- \theta(\gamma)) |^2 dx } \\
&=& \| a \|_1 \| f \|_2
\eean
Now if $(S_n)_{n \in \BN}$ is a sequence of finite subsets of $\BZ^d$ for which $S_n \subset S_{n+1}$ and $\cup_{n \in \BN} S_n = \BZ^d$, then 
$$
(f \circ_\theta a)(x) = \lim_{n \to \infty} \sum_{\gamma \in S_n} a(\gamma) f(x- \theta(\gamma)).
$$
This series is absolutely bounded by the bound $\|a \|_1 \|f \|_2$ and it therefore follows that $\|f \circ_\theta a \|_2 \leq \|a \|_1 \| f \|_2 < \infty$ and the sum in equation (\ref{rdbracmaeq-2}) converges unconditionally in $L^2(\BR^d)$.

\item This is a routine calculation which is done in Lemma 2.1.6 of \cite{w1}. Let us verify (\ref{eqrmod2}) to see how these calculations are done. For $x \in \BR^d$,
\bean
((f \circ_\theta a) \circ_\theta b)(x) &=& \sum_{\gamma \in \BZ^d} (f \circ_\theta a)(x - \theta(\gamma))b(\gamma) \\
&=& \sum_{\gamma \in \BZ^d} \sum_{\gamma' \in \BZ^d} f(x - \theta(\gamma + \gamma')) a(\gamma') b(\gamma) \\
&=& \sum_{\beta \in \BZ^d} f(x - \theta(\beta)) \sum_{\gamma \in \BZ^d} a(\beta - \gamma) b(\gamma) \\
&=& (f \circ_\theta (a * b))(x).
\eean
verifying (\ref{eqrmod2}).

\item To prove part 4, one first uses the Cauchy-Schwarz inequality to verify that $[f,g]_\theta \in l^\infty(\BZ^d)$. It then follows that if $f_n \to f$ and $g_n \to g$ are convergent sequences in $L^2(\BR^d)$, then $[f_n,g_n]_\theta \to [f,g]_\theta$ in $l^\infty(\BZ^d)$. Because $C_0(\BZ^d)$ is the completion of $C_c(\BZ^d)$ with respect to the supremum norm, it follows that $[f,g]_\theta \in C_0(\BZ^d)$.

\item This is a routine calculation which is done in Lemma 2.1.14 of \cite{w1}.
\end{enumerate}
\hfill \end{Proof}

\begin{Lemma} \label{rdftbraclma2-2}
Suppose that $p,q \in L^2(\widehat{\BR^d})$ and $\theta  : \BZ^d \to \BR^d$ satisfies $\theta(\gamma) = A_\theta \iota \gamma$, for a linear transformation 
$A_\theta: \BR^d \to \BR^d$. Let $\hl p,q \hr_\theta$ be the function on $\BT^d$ given by
\be \label{eq:hlpqhr}
\hl p,q \hr_\theta (\zeta) := \frac{1}{\det (A_\theta)} \sum_{\hat{\theta} (\xi) = \zeta} \overline{p(\xi)} q(\xi)
\ee
for almost every $\zeta \in \BT^d$. Then the above sum converges absolutely almost everywhere and $ \hl p,q \hr_{\theta} \in l^1(\BT^d)$. We furthermore have that $\{ \xi \in \widehat{\BR^d} : \hat{\theta}(\xi) = 0 \} = (A_\theta^*)^{-1} \BZ^d = \mathrm{Ann }\theta(\BZ^d)$. We hence can write
\bea
\label{eq:hlpqhr2} \hl p,q \hr_\theta (\hat{\theta} (\xi)) 
                     &=& \frac{1}{\det (A_\theta)} \sum_{\beta \in \mathrm{Ann} \ \theta(\BZ^d)} \overline{p(\xi + \beta)} q(\xi + \beta) \\
  \label{eq:hlpqhr3} &=& \frac{1}{\det (A_\theta)} \sum_{\beta \in (A_\theta^*)^{-1} \BZ^d} \overline{p(\xi + \beta)} q(\xi + \beta)
\eea
for almost every $\xi \in \BR^d$. 
\end{Lemma}
\begin{Proof}{}
We shall first consider the function $\bar{p}q$ which is contained in $L^1(\widehat{\BR^d})$. Define a group isomorphism $I : \BZ^d \times \BT^d \to \BR^d$ by setting $I (\beta, \zeta) = \beta + \zeta$ for $\beta \in \BZ^d$, $\zeta \in \BT^d$. The group isomorphism $I$ extends to a measure space isomorphism from $\BZ^d \times \BT^d$ equipped with the product measure $ \mu_{\BZ^d} \times \mu_{\BT^d} $ onto the measure space $(\BR^d, \mu_{\BR^d})$. Here $\mu_{\BZ^d}$ is counting measure on $\BZ^d$, $\mu_{\BT^d}$ is Lebesgue measure normalised so that $\mu_{\BT^d}(\BT^d) = 1$, and $\mu_{\BR^d}$ is Lebesgue measure on $\BR^d$.

Since $\bar{p}q \in L^1(\widehat{\BR^d})$, it follows that
$$
p \circ I \in L^1(\BZ^d \times \BT^d,  \mu_{\BZ^d} \times \mu_{\BT^d}).
$$
It follows directly from Fubini's Theorem that the series
$$
\sum_{\beta \in \BZ^d} \overline{p(\xi + \beta)} q(\xi + \beta)
$$
is absolutely convergent for almost every $\xi \in \widehat{\BR^d}$. Therefore the series
\be \label{sumhlpqhr}
\sum_{\beta \in \mathrm{Ann }\theta(\Gamma)} \overline{p(\xi + \beta)} q(\xi + \beta)
\ee
is also absolutely convergent for almost every $\xi \in \widehat{\BR^d}$.

For any $\beta \in \mathrm{Ann} \ \theta(\BZ^d)$ and $\xi \in \widehat{\BR^d}$, we have that $\hat{\theta}(\xi + \beta) = \hat{\theta}(\xi)\hat{\theta}(\beta) = \hat{\theta}(\xi)$. We therefore have that $\{ \xi \in \widehat{\BR^d} : \hat{\theta}(\xi) = 0 \} = \mathrm{Ann }\theta(\BZ^d)$, verifying equation (\ref{eq:hlpqhr2}). Because (\ref{sumhlpqhr}) converges absolutely, it follows from (\ref{eq:hlpqhr2}) that $ \hl p,q \hr_\theta (\zeta)$ converges absolutely for almost every $\zeta \in \BT^d$. We know from Lemma \ref{thetalma} that $\mathrm{Ann} \ \theta(\BZ^d) = (A_\theta^*)^{-1} \BZ^d$, and this gives us the equality between equations (\ref{eq:hlpqhr2}) and (\ref{eq:hlpqhr3}).
\end{Proof}

\begin{Definition} \label{rdftbpma}
Suppose that $p,q \in L^2(\widehat{\BR^d})$ and $\theta  : \BZ^d \to \BR^d$ satisfies $\theta(\gamma) = A_\theta \iota \gamma$, for a linear transformation 
$A_\theta: \BR^d \to \BR^d$. We call $\hl p,q \hr_\theta$ the {\em Fourier transformed bracket product} associated with $\theta$, where $\hl p,q \hr_\theta$ 
is defined by equation (\ref{eq:hlpqhr}).
Suppose $b \in C(\BT^d)$, we define the {\em Fourier transformed module action} associated with $\theta$ to be the function on $\widehat{\BR^d}$ denoted by $p 
\widehat{\circ_\theta} b$ and given by
\be
(p \widehat{\circ_\theta} b)(\xi) := p(\xi) b( \hat{\theta} (\xi))
\ee
for almost every $\xi \in \widehat{\BR^d}$.
\end{Definition}

We remark that because $b$ is contained in $C(\BT^d)$, it follows that $\sup_{\xi \in \BR^d} b(\hat{\theta} \xi) < \infty$. Because $| (p \widehat{\circ_\theta} b)(\xi) | \leq | p(\xi) \sup_{\xi' \in \BR^d} b( \hat{\theta}(\xi'))|$, it follows that $ p \widehat{\circ_\theta} b \in L^2(\BR^d)$. It is the case that $\widehat{\circ_\theta}$ makes $L^2(\widehat{\BR^d})$ into a right $C(\BT^d)$-module because $\widehat{\circ_\theta}$ consists of pointwise multiplication.

The proof to Lemma \ref{rdftbracident} is a routine calculation which is similar to part 5 of Proposition \ref{bptheorem} (see \cite[Lemma 2.1.15]{w1}).

\begin{Lemma}
\label{rdftbracident}
Suppose that $p,q, r \in L^2(\widehat{\BR^d})$, $b \in C(\BT^d)$, $ \alpha, \beta \in \BC$, and $\theta$ is an embedding of $\BZ^d$ in $\BR^d$, then
\begin{enumerate}
\item ${\hl p, \alpha q + \beta r \hr}_{\theta} = \alpha \hl p,q \hr_{\theta}  + \beta \hl p,r \hr_{\theta} $;
\item ${\hl  p, q \widehat{\circ_\theta} b  \hr}_{\theta} = \hl p,q \hr_{\theta} b$, where $b$ is a function on $\BT^d$ for which $q \widehat{\circ_\theta} b \in L^2(\BR^d)$ whenever $q \in L^2(\BR^d)$ (this is the case when $b \in C(\BT^d)$ or $b \in L^\infty(\BT^d)$);
\item ${\hl  p,q \hr}_{\theta}^* = \hl q,p \hr_{\theta}$, where the unary operation $\cdot^*$ is complex conjugation on $C(\BT^d)$.
\end{enumerate}
\end{Lemma}

The following result justifies our use of the term ``Fourier transformed Bracket product''.

\begin{Lemma}
\label{lma:rdfouhilb}
For $f,g \in L^2(\BR^d)$ and $a \in l^1(\BZ^d)$, suppose that $[f,g]_\theta \in l^1(\BZ^d)$, then
\bea
\label{eq:rdfh1} \hat{f} \widehat{\circ_\theta} \hat{a} &=& \mathcal{F}_{\BR^d} (f \circ_\theta a) \\
\label{eq:rdfh2} \mbox{and } \ \ \hl \hat{f}, \hat{g} \hr_{\theta} &=& \mathcal{F}_{\BZ^d} ( [f,g]_\theta )
\eea
where $\mathcal{F}_{\BR^d}$ is the Fourier transform on $\BR^d$, and $\mathcal{F}_{\BZ^d}$ is the Fourier transform on $\BZ^d$. Furthermore, if $p,q \in L^2(\widehat{\BR^d})$, then for $\gamma \in \BZ^d$
\be 
\label{eq:rdfh3} \mathcal{F}_{\BT^d}(\hl p,q \hr_{\hat{\theta}})(\gamma) = [\check{p}, \check{q}]_\theta (\gamma)
\ee
where $\check{p}, \check{q} \in L^2(\BR^d)$ are the inverse Fourier transforms of $p,q$ and $\mathcal{F}_{\BT^d}$ is the Fourier transform on $\BT^d$. 
\end{Lemma}
\begin{Proof}{}
We know from equation (\ref{rdbpmarem2}) that for $\gamma \in \BZ^d$, $(\check{p}^* * \check{q})(\theta(\gamma)) = [\check{p},\check{q}]_\theta (\gamma)$, so we make use of Lemma \ref{thetalma} to calculate
\bean
(\mathcal{F}_{\BT^d} \hl p, q \hr_{\theta})(\gamma)
&=& \int_{\zeta \in \BT^d}  \hl p, q \hr_{\theta} (\zeta) (-\gamma, \zeta) d \zeta \\
&=& \int_{\zeta \in \BT^d} \frac{1}{\det (A_\theta)} \sum_{\hat{\theta} \xi = \zeta} \overline{p(\xi)} q(\xi) (-\gamma, \zeta) d \zeta \\
&=& \int_{\zeta \in \BT^d} \frac{1}{\det (A_\theta)} \sum_{\hat{\iota} A_\theta^* \xi = \zeta} \overline{p(\xi)} q(\xi) (-\gamma, \hat{\iota} A_\theta^* \xi) 
d \zeta \\
\mbox{now set $\eta := A_\theta^* \xi$},&\ & \\
&=& \int_{\zeta \in \BT^d} \frac{1}{\det (A_\theta)} \sum_{\hat{\iota} \eta = \zeta} \overline{p( (A_\theta^*)^{-1} \eta)} q( (A_\theta^*)^{-1} \eta)  
(-\gamma, \hat{\iota} \eta) d \zeta \\
&=& \int_{\eta \in \BR^d} \frac{1}{\det (A_\theta)} \overline{p((A_\theta^*)^{-1} \eta)} q((A_\theta^*)^{-1} \eta)  (-\gamma, \hat{\iota} \eta) d \eta \\
&=& \int_{\xi \in \BR^d} \overline{p(\xi)} q(\xi)  (-\gamma, \hat{\iota} A_\theta^* \xi) d \xi \\
&=& \int_{\xi \in \BR^d} \overline{p(\xi)} q(\xi)  (\theta (\gamma), \xi) d \xi \\
&=& (\check{p}^* * \check{q})(\theta (\gamma)) \\
&=& [\check{p},\check{q}]_\theta (\gamma).
\eean
Now because we assumed that $[f,g]_\theta \in l^1(\BZ^d)$, the Fourier transform of $[f,g]_n$ is defined. Therefore $ \hl \hat{f}, \hat{g} \hr_\theta = \mathcal{F}_{\BZ^d} ( [f,g]_\theta )$.

Now we also have
\bean
(\mathcal{F}_{\BR^d}(f \circ_\theta a))(\xi) &=& \int_{\BR^d} \sum_{\gamma \in \BZ^d} a(\gamma) f(x- \theta( \gamma)) (x, \xi) dx \\
&=& \int_{\BR^d} \sum_{\gamma \in \BZ^d} a(\gamma) f(x) (x + \theta( \gamma), \xi) dx \\
&=& \sum_{\gamma \in \BZ^d} a(\gamma) (\theta( \gamma), \xi) \int_{\BR^d} f(x) (x, \xi) dx \\
&=& \sum_{\gamma \in \BZ^d} a(\gamma) (\gamma, \hat{\theta} (\xi)) \int_{\BR^d} f(x) (x, \xi) dx \\
&=& \hat{a}(\hat{\theta}(\xi)) \hat{f}(\xi) \\
&=& (\hat{f} \widehat{\circ_\theta} \hat{a})(\xi).
\eean
verifying the result.
\end{Proof}

\begin{Proposition}
\label{thetainnerprodmod}
Let $\theta : \BZ^d \to \BR^d$ be an embedding. The space $C_c(\BR^d)$ is an inner product $C_c(\BZ^d)$-module with operations $[ \;, \: ]_\theta$ and $\circ_\theta$. For $f \in C_c(\BR^d)$ and an embedding $\theta:\BZ^d \to \BR^d$, let
\be \label{eq:hmnx1}
\|f\|_{X_\theta} := \sup_{\zeta \in \BT^d} \sqrt{ \widehat{[f,f]_\theta} (\zeta) }.
\ee
Then $\| \cdot \|_{X_\theta}$ is a norm equal to the Hilbert module norm on $C_c(\BR^d)$.
\end{Proposition}
\begin{Proof}{}
In part 3 of Proposition \ref{bptheorem}, we showed that $C_c(\BR^d)$ is a right $C_c(\BZ^d)$ module with the operation $\circ_\theta$.

From part 5 of Proposition \ref{bptheorem}, properties 1, 2 and 3 of Definition \ref{def:hmodule} are satisfied for $C_c(\BR^d)$ with the operations $\circ_\theta$ and $[ \; , \:]_\theta$. Suppose that $f \in C_c(\BR^d)$, then by Lemma \ref{lma:rdfouhilb},
$$
(\mathcal{F}_{\BZ^d} [f,f]_\theta )(\zeta) = \hl \hat{f}, \hat{f} \hr_{\theta} (\zeta) = \frac{1}{\det (A_\theta)} \sum_{\beta \in (A_\theta^*)^{-1} \BZ^d} 
\overline{\hat{f}(\zeta + \beta)} \hat{f}(\zeta + \beta)
$$
which is non-negative for all $\zeta \in \BT^d$. So $\mathcal{F}_{\BZ^d} [f,f]_\theta$ is a positive element of the $C^*$-algebra $C(\BT^d)$. It is well known (see for example \cite[Proposition VII.1.1]{dv1}) that $C^*(\BZ^d)$ and $C(\BT^d)$ are isomorphic and that this isomorphism is given by the Fourier transform on the subalgebra $l^1(\BZ^d)$ which is dense in $C^*(\BZ^d)$. This implies that $[f,f]_\theta$ is a positive element of $C^*(\BZ^d)$, verifying Property 4 of Definition \ref{def:hmodule}. Now if $[f,f]_\theta = 0$, then for all $\gamma \in \BZ^d$, $\int_{\BR^d} \overline{f(x- \gamma)} f(x) dx = 0$, so in particular $\int_{\BR^d} \overline{f(x)} f(x) dx =0$, so $\|f \|_2 = 0$, implying that $f = 0$. This verifies Property 5 of Definition \ref{def:hmodule}. It therefore follows that $C_c(\BR^d)$ is a right inner product $C_c(\BZ^d)$-module. Because we know that $C_c(\BR^d)$ is an inner-product $C_c(\BZ^d)$-module, the fact that $\| \cdot \|_{X_\theta}$ is a norm directly follows from \cite[Corollary 2.7]{rw1}.It is an immediate consequence of Definition \ref{def:hmodule} and the definition of $\| \cdot \|_{X_\theta}$ that $\| \cdot \|_{X_\theta}$ is the Hilbert module norm.
\end{Proof}

It is worth noting that Proposition \ref{thetainnerprodmod} is a special case of one of the main results in \cite{r1} (see also \cite{rw1}). The results in these references are more general in that $\BZ^d$ and $\BR^d$ are replaced by a closed subgroup $H$ of a locally compact group $G$. The fact that $\BZ^d$ is abelian enables us to use the Fourier transform in the proof for Proposition \ref{thetainnerprodmod} to provide a simpler proof of the positivity of the $C_c(\BZ^d)$-valued inner product.

\begin{Lemma}
\label{lma:hmnormx}
The norm $\| \cdot \|_{X_\theta}$ can be expressed as
\bea 
\|f\|_{X_\theta} &=& \sup_{\zeta \in \BT^d} \sqrt{ \hl \hat{f}, \hat{f} \hr_{\theta} (\zeta)} \\    \label{eq:hmnx2}
&=& \sup_{\zeta \in \BT^d} \sqrt{ \frac{1}{\det (A_\theta)} \sum_{\beta \in (A_\theta^*)^{-1} \BZ^d} \overline{\hat{f}(\zeta + \beta)} \hat{f} (\zeta + \beta) }.
\eea
Furthermore, for all $f \in C_c(\BR^d)$, $\|f \|_2 \leq \|f \|_{X_\theta}$, implying that the completion of $C_c(\BR^d)$ with respect to the norm $\| \cdot \|_{X_\theta}$ is contained in $L^2(\BR^d)$.
\end{Lemma}
\begin{Proof}{}
From equation (\ref{eq:rdfh2}) of Lemma \ref{lma:rdfouhilb}, $\widehat{[f,f]_\theta} = \hl \hat{f}, \hat{f} \hr_\theta$, so
\bean
\|f\|_{X_\theta}^2 &=& \sup_{\zeta \in \BT^d} \hl \hat{f}, \hat{f} \hr_{\theta} (\zeta) \\
&=& \sup_{\zeta \in \BT^d} \frac{1}{\det (A_\theta)} \sum_{\beta \in (A_\theta^*)^{-1} \BZ^d} \overline{\hat{f}(\zeta + \beta)} \hat{f} (\zeta + \beta),
\eean
proving equation (\ref{eq:hmnx2}).

Now we have for $f \in C_c(\BR^d)$,
\bean
\| f \|_2^2 &=& \int_{\BR^d} f(x) \overline{f(x)} dx \\
&=& \int_{\BR^d} \hat{f}(x) \overline{\hat{f}(x)} dx  \\
&=& \int_{\BT^d} \frac{1}{\det (A_\theta)} \sum_{\beta \in (A_\theta^*)^{-1} \BZ^d} \overline{\hat{f}(\zeta + \beta)} \hat{f} (\zeta + \beta) d \zeta  \\
&=& \int_{\BT^d} \hl \hat{f}, \hat{f} \hr_\theta (\zeta) d \zeta \\
&\leq& \sup_{\zeta \in \BT^d} \hl \hat{f}, \hat{f} \hr_\theta (\zeta) \\
&=& \| f \|_{X_\theta}^2 \\
\eean
proving the assertion.
\end{Proof}


We shall now define the Hilbert modules $X_\theta$ and $\hat{X_\theta}$, we will then show that they are Hilbert modules.

\begin{Definition} \label{rddefx}
Let $\theta : \BZ^d \to \BR^d$ be an embedding. We define the {\em bracket product Hilbert $C^*(\BZ^d)$-module} $X_\theta$ as follows: Let the 
linear space $X_\theta$ be the completion of $C_c(\BR^d)$ with respect to the norm $\| \cdot \|_{X_\theta}$ that was defined in Lemma \ref{lma:hmnormx}. We 
equip $X_\theta$ with the bracket product ${[\; ,\:]}_\theta : X_\theta \times X_\theta \to C^*(\BZ^d)$ and module action $\circ_\theta : X_\theta \times 
C^*(\BZ^d) \to X_\theta$ as defined in Definition \ref{rdbpma}.

We define the {\em Fourier transformed bracket product Hilbert $C(\BT^d)$-module} $\hat{X_\theta}$ as follows: Let the linear space $\hat{X_\theta}$ be the 
image of $X_\theta$ under the Fourier transform on $\BR^d$. We equip $\hat{X_\theta}$ with the Fourier transformed bracket product ${\hl \; ,\: \hr}_{\theta} 
: \hat{X_\theta} \times \hat{X_\theta} \to C(\BT^d)$ and module action $\widehat{\circ_\theta} : \hat{X_\theta} \times C(\BT^d) \to \hat{X_\theta}$, as 
defined in Definition \ref{rdftbpma}.
\end{Definition}

We know from Lemma \ref{lma:hmnormx} that we can continuously embed $X_\theta$ in $L^2(\BR^d)$, this means that it makes sense to talk about the image of $X_\theta$ under the Fourier transform on $\BR^d$. In Definition \ref{rdbpma} we defined the bracket product on $L^2(\BR^d)$, so we can use this to describe the bracket product on $X_\theta$. It is very useful that we have also constructed a Hilbert module $\hat{X_\theta}$ in the Fourier domain, because in practice we can use the Gelfand transform to represent an element of $C^*(\BZ^d)$ as a continuous function on $\BT^d$, and use the Fourier transformed module action $\widehat{\circ_\theta}$.

Theorem \ref{rdhmexists} is the main theorem of this section and states that Definition \ref{rddefx} defines a Hilbert module. The main idea of the proof is to use the completion process described in Lemma \ref{comphm} to obtain a Hilbert $C^*(\BZ^d)$-module from the inner product $C_c(\BZ^d)$-module $C_c(\BR^d)$.

\begin{Theorem}
\label{rdhmexists}
Let $\theta : \BZ^d \to \BR^d$ be an embedding. Let $X_\theta \subset L^2(\BR^d)$ be the completion of $C_c(\BR^d)$ with respect to the norm $\| \cdot \|_{X_\theta}$ as defined in Definition \ref{rddefx}. Let $\hat{X_\theta} \subset L^2(\widehat{\BR^d})$ be the image of $X_\theta$ under the Fourier transform as defined in Definition \ref{rddefx}. We have that
\begin{enumerate}
\item The space $X_\theta$ is a full Hilbert $C^*(\BZ^d)$-module with the inner product $[ \; , \:]_\theta$. The module action of $C^*(\BZ^d)$ on $X_\theta$ has the property that for $f \in X_\theta$ and $a \in l^1(\BZ^d)$ (note that $l^1(\BZ^d) \subset C^*(\BZ^d)$), the module action of $a$ on $f$ is given by the module action associated with the bracket product $f \circ_\theta a$ as defined in Definition \ref{rddefx}. We furthermore have that if $(f_n)_{n \in \BZ}$ is a convergent sequence in $X_\theta$, then it is a convergent sequence in $L^2(\BR^d)$. 
\item The operations $\hl \; , \: \hr_{\hat{\theta}}$ and $\widehat{\circ_\theta}$ make $\hat{X_\theta}$ into a full Hilbert $C(\BT^d)$-module;
\item If we identify the isomorphic $C^*$-algebras $C^*(\BZ^d)$ and $C(\BT^d)$, then the Fourier transform on $\BR^d$ is a unitary isomorphism between $X_\theta$ and $\hat{X_\theta}$. 
\end{enumerate}
\end{Theorem}
\begin{Proof}{}
\begin{enumerate}
\item It follows from Lemma \ref{comphm} and Proposition \ref{thetainnerprodmod} that the space $X_\theta$ is a Hilbert $C^*(\BZ^d)$-module with the inner product $[ \; , \:]_\theta$. It also follows from Lemma \ref{comphm} that the module action of $C^*(\BZ^d)$ on $X_\theta$ has the property that for $f \in X_\theta$ and $a \in l^1(\BZ^d)$, the module action of $a$ on $f$ is given by the module action associated with the bracket product $f \circ_\theta a$. We know from Lemma \ref{lma:hmnormx} that if $f \in C_c(\BR^d)$, then $\| f \|_2^2 \leq \| f \|_{X_\theta}^2$. It follows that if $f, g \in C_c(\BR^d)$ and $\varepsilon > 0$, then $\|f - g \|_{X_\theta} < \varepsilon$ implies $\| f - g \|_2 < \varepsilon$. Therefore if $(f_n)_{n \in \BN}$ is a convergent sequence in $X_\theta$, then it is a convergent sequence in $L^2(\BR^d)$.

To see that $X_\theta$ is full, choose an element $\phi \in X_\theta$ for which $[ \phi, \phi ]_\theta = \mathbf{1}$. 
A possible choice for $\phi$ would be any continuous compactly supported wavelet or scaling function corresponding to the translations $\theta(\BZ^d)$. For example, let $\varphi$ be the Daubechies scaling function (see \cite{daub1}). The Daubechies scaling function $\varphi$ is a continuous compactly supported function on $\BR$ satisfying $[ \varphi, \varphi ]_\theta = \mathbf{1}$. We can define a continuous compactly supported function $\phi$ on $\BR^d$ or which $[ \phi, \phi ]_\theta = \mathbf{1}$ as follows. We let
$$
\phi(\theta(x_1, x_2, \ldots, x_d)) = \varphi(x_1) \varphi(x_2) \ldots \varphi(x_d).
$$
Then $[ \phi, \phi ]_\theta = \mathbf{1}$. Therefore for all $a \in C^*(\BZ^d)$, it is the case that
\bean
{ [ \phi , \phi \circ_\theta a ]_\theta} &=& [\phi, \phi]_\theta \circ_\theta a \\
&=& a.
\eean
Because $\phi \circ_\theta a \in X_\theta$, it follows that $C^*(\BZ^d) = [  X_\theta , X_\theta ]_\theta$, verifying that $X_\theta$ is full.

\item Let $\widehat{C_c(\BR^d)}$ be the image of $C_c(\BR^d)$ under the Fourier transform. Let $\widehat{C_c(\BZ^d)} \subset C(\BT^d)$ be the image of $C_c(\BZ^d)$ under the Fourier transform. By Lemma \ref{lma:rdfouhilb}, if $p,q \in \widehat{C_c(\BR^d)}$, then $\hl p,q \hr_{\theta} \in \widehat{C_c(\BZ^d)}$. It also follows from Lemma \ref{lma:rdfouhilb} that $\widehat{C_c(\BR^d)}$ is a right $\widehat{C_c(\BZ^d)}$-module with module action 
$\widehat{\circ_\theta}$.

From Lemma \ref{rdftbracident}, Properties 1, 2 and 3 of Definition \ref{def:hmodule} are satisfied for $\widehat{C_c(\BR^d)}$ with the operations $\widehat{\circ_\theta}$ and $\hl \; , \: \hr_{\theta}$. Properties 4 and 5 of Definition \ref{def:hmodule} are satisfied by the same argument as for $C_c(\BR^d)$. Therefore the operations $\hl \; , \: \hr_{\theta}$ and $\widehat{\circ_\theta}$ make $\widehat{C_c(\BR^d)}$ into an inner product $\widehat{C_c(\BZ^d)}$-module. The completion of $\widehat{C_c(\BR^d)}$ with respect to the Hilbert module norm is equal to $\hat{X_\theta}$. So by Lemma \ref{comphm}, the operations $\hl \; , \: \hr_{\theta}$ and $\widehat{\circ_\theta}$ make $\hat{X_\theta}$ into a Hilbert $C(\BT^d)$-module, and if $p,q \in \hat{X_\theta}$, then $\hl p,q \hr_{\theta} \in C(\BT^d)$. The fact that $\hat{X_\theta}$ is full follows from the fact that $X_\theta$ is full.

\item By Definition \ref{def:isohm}, we need to show that the Fourier transform $\mathcal{F}_{\BR^d}$ is a unitary adjointable operator from $X_\theta$ to $\hat{X_\theta}$. We want to show that $\mathcal{F}$ is a unitary element of $\mathcal{L}(X_\theta,\hat{X_\theta})$. From Lemma \ref{lma:rdfouhilb}, if $f \in X_\theta$ and $p \in \hat{X_\theta}$, then
\be
{\hl \mathcal{F}_{\BR^d} ( f ) , p \hr}_{\theta} = \mathcal{F}_{\BZ^d} \left( [ f , \mathcal{F}_{\BR^d}^* (p)]_\theta \right) .
\ee

We therefore have that $\mathcal{F}_{\BR^d}$ is adjointable with Hilbert module adjoint $\mathcal{F}_{\BR^d}^*$. We have that $\mathcal{F}_{\BR^d} \mathcal{F}_{\BR^d}^* = \mathbf{1}_X$ and $\mathcal{F}_{\BR^d}^* \mathcal{F}_{\BR^d} = \mathbf{1}_{\hat{X}}$, so $\mathcal{F}_{\BR^d}$ is unitary.

\end{enumerate}
\hfill \end{Proof}

Some of the theory described above can be generalised considerably. Consider a multiresolution structure $(\Gamma, \dil)$ with index $m$ acting on the Hilbert space $L^2(G)$, for $G$ a locally compact group. Assume that $\Gamma$ is a subgroup of $G$. Suppose that $\theta : \Gamma \to G$ is an embedding (we know that such a $\theta$ exists because $\Gamma$ is assumed to be a subgroup of $G$). Then $\theta(\Gamma)$ is also a subgroup of $G$. Associated with $\theta$ we define a representation $\pi^\theta$ of $\Gamma$ on $L^2(G)$ by
$$
(\pi^\theta_\gamma f)(x) = f(x \theta(\gamma^{-1}) ), \ \ \mbox{for $x \in G$, $\gamma \in \Gamma$, $f \in L^2(G)$.}
$$
Let $\Delta$ be the modular function on $G$ and let $\delta$ be the modular function on $\Gamma$. For $f,g \in C_c(G)$, $a \in C_c(\Gamma)$, define
\bea
(f \circ_\theta a)(s) := \int_\Gamma f(s \theta( t^{-1})) \delta(t^{-1})  a(t) d \mu_\Gamma (t) ; \\
{[f,g]_\theta}(t) := \sqrt{ \frac{\Delta(t)}{\delta(t)}} \int_G \overline{f(r)} g(r \theta(t)) d \mu_G (r).
\eea
Then from \cite[Theorem C.23]{rw1}, the completion of $C_c(G)$ with respect to the Hilbert module norm 
\be
\| f \|_X = \| [f,f]_\theta \|_*^\half, \ \ \ f \in C_c(G);
\ee
is a Hilbert $C^*(\Gamma)$-module.

At this level of generality there is no longer an obvious way of defining the Fourier transform. There are impediments to proving an analogue of Lemma \ref{rdftbraclma2-2} because it identifies $\BT^d$ with the cube $[-1/2,1/2)^d$, and in the more general situation it is difficult to come up with an analogue of the cube which has all of the required properties (for example $[-1/2,1/2)^d$ is a measurable subset of $\BR^d$ for which a neighbourhood of $[-1/2,1/2)^d$ contains $0$, and we can tile $\BR^d$ with translations of $[-1/2,1/2)^d$ by elements of $\BZ^d$). 

\section{Incorporating the dilation} \label{s4}

In this section we shall relate the Hilbert $C^*(\BZ^d)$-modules defined in the previous section to wavelets. Let $(\Gamma=\BZ^d,\dil)$ be the standard multiresolution structure associated with a dilation matrix $\tilde{\dil}$. 

Let us first examine the Harmonic analysis of the groups $\Gamma^n$, (recall from Definition \ref{def:mrs} that $\Gamma^n = \dil^n \BZ^d \dil^{-n}$). For each integer $n$, we can define an injective group homomorphism $\iota_n: \BZ^d \to \BR^d$ by
\be \label{eq:iotainrd}
\iota_n(\gamma) = \tilde{\dil}^{-n} \gamma
\ee
for $\gamma \in \BZ^d$. Let $\iota$ be the natural embedding of $\BZ^d$ into $\BR^d$ and note that $\iota_n(\gamma) = \tilde{\dil}^{-n} \iota \gamma$, and $\iota_0 = \iota$. The image of $\iota_n$ is the group $\tilde{\dil}^{-n} \BZ^d$ and $\tilde{\dil}^{-n} \BZ^d$ acts by translations on $L^2(\BR^d)$ in the same way as $\Gamma^n$. We regard $\Gamma^n$ as a subgroup of $\BR^d$ by identifying it with $\tilde{\dil}^{-n} \BZ^d$. Define another homomorphism $\hat{\iota_n}: \widehat{\BR^d} \to \BT^d$ by
$$
\hat{\iota_n} \xi = \hat{\iota} \tilde{\dil}^{*-n} \xi
$$
where $\xi \in \widehat{\BR^d}$, and $\tilde{\dil}^*$ is the adjoint of $\tilde{\dil}$. A calculation using basic properties of characters on locally compact abelian groups (\cite{w1}, p34) verifies that for $\gamma \in \BZ^d$, $\xi \in \widehat{\BR^d}$, $(\iota_n \gamma, \xi) = (\gamma, \hat{\iota_n} \xi)$. Let us now examine the annihilator of $\Gamma^n$. Using the definition of an annihilator it can be shown using a routine calculation (\cite{w1}, p35) that $\mathrm{Ann }\Gamma^n = \tilde{\dil}^{*n} \BZ^d. $

We furthermore have that $\xi \in \mathrm{Ann }\Gamma^n$ if and only if for all $\gamma \in \BZ^d$, $( \gamma, \hat{\iota_n} \xi) = 1$. This means that if $a$ is a function on $\BT^d$, then $a \circ \hat{\iota_n}$ is an $\mathrm{Ann }\Gamma^n$-periodic function on $\widehat{\BR^d}$.

In Theorem \ref{rdhmexists} we have shown how to construct a Hilbert module $X_\theta$ from an embedding $\theta: \BZ^d \to \BR^d$. Recall from Lemma \ref{thetalma} that associated with this embedding is a linear transformation $A_\theta : \BR^d \to \BR^d$ for which $\theta = A_\theta \iota$, where $\iota$ is the natural embedding of $\BZ^d$ in $\BR^d$. We also showed in Lemma \ref{thetalma} that there is a dual homomorphism $\hat{\theta}: \widehat{\BR^d} \to \BT^d$ given by $\hat{\theta} = \hat{\iota} A_\theta^*$. We have that $A_{\iota_n} = \tilde{\dil}^{-n}$.

\begin{Definition} \label{def:noperations}
Suppose that $\tilde{\dil}$ is a dilation matrix, and $n$ is an integer.
For $f,g \in L^2(\BR^d)$, $a \in l^1(\BZ^d)$, $\gamma \in \BZ^d$, and $x \in \BR^d$, define the {\em $n$th level bracket product} $[ \; , \: ]_n$ and {\em $n$th level module action} $\circ_n$ to be
\bea
{[f,g]_n(\gamma)} &:=& [f,g]_{\iota_n}(\gamma) = \int_{\BR^d} \overline{f(x - \tilde{\dil}^{-n} \iota( \gamma))} g(x) dx, \\
f \circ_n a &:=& f \circ_{\iota_n} a = \sum_{\gamma \in \BZ^d} a(\gamma) \pi^n_\gamma f. 
\eea
For $p,q \in L^2(\widehat{\BR^d})$ and $b \in C(\BT^d)$, we define the {$n$th level Fourier transformed bracket product} $\hl p, q \hr_n$ and {\em $n$th level Fourier transformed module action} $\widehat{\circ_n}$ to be
\bea
\hl p, q \hr_n &:=& \hl p, q \hr_{\iota_n}, \\
p \widehat{\circ_n} b &:=& p \widehat{\circ_{\iota_n}} b.
\eea
\end{Definition}

\begin{Lemma} \label{pqeqlma}
For $p,q \in \widehat{\BR^d}$, $\zeta \in \BT^d$, if we identify $\BT^d$ with the cube $[ -\half , \half)^d$, then
\be \label{ftbpuseful}
\hl p, q \hr_n (\zeta) = m^n \sum_{\beta \in (\tilde{\dil}^*)^n \BZ^d} \overline{ p \left( (\tilde{\dil}^*)^n \zeta + \beta \right) } q \left( (\tilde{\dil}^*)^n \zeta + \beta \right) .
\ee
\end{Lemma}
\begin{Proof}{}
Because $A_{\iota_n} = \tilde{\dil}^{-n}$, it follows that $\det(A_{\iota_n}) = m^{-n}$ and that $(A_{\iota_n}^*)^{-1} \BZ^d = (\tilde{\dil}^*)^n \BZ^d$. Now from Lemma \ref{rdftbraclma2-2}, $(\tilde{\dil}^*)^n \BZ^d = \{ \xi \in \widehat{\BR^d} : \hat{\theta}(\xi) = 0 \}$. It is therefore the case that $\{ \xi : \hat{\iota_n}(\xi) = \zeta \} = \{ (\tilde{\dil}^*)^n \zeta + \beta : \beta \in (\tilde{\dil}^*)^n \BZ^d \}$. The result now follows from equation (\ref{eq:hlpqhr}).
\end{Proof} 

The following Lemma tells us how the operations we have just defined relate to each other when $n$ changes.
\begin{Lemma}
\label{dilcalc}
Suppose that $n$ is an integer. We have
\begin{enumerate}
\item If $f,g \in L^2(\BR^d)$, then
\be \label{dilident}
{[f,g]_n} = [ \dil^{-n} f, \dil^{-n} g]_0.
\ee
\item If $p,q \in L^2(\widehat{\BR^d})$, then
\be \label{dilhatident}
\hl p, q \hr_n = \hl \hat{\dil}^n p, \hat{\dil}^n q \hr_0.
\ee
\item If $f \in L^2(\BR^d)$, $a \in l^1(\BZ^d)$, then 
\be \label{dilmaident}
\dil^{-n} (f \circ_n a) = (\dil^{-n} f) \circ_0 a.
\ee
\item If $p \in L^2(\widehat{\BR^d})$, $b \in C(\BT^d)$, then
\be \label{dilhatmaident}
\hat{\dil}^{-n} (p \widehat{\circ_n} b) = (\hat{\dil}^{-n} p) \widehat{\circ_0} b).
\ee
\end{enumerate}
\end{Lemma}
\begin{Proof}{}
The calculations that verify this result involve substituting (\ref{dildef}) and (\ref{dilhatdef}) into the relevant equations in Definition \ref{def:noperations} and Lemma \ref{pqeqlma}. The interested reader is referred to \cite[Lemma 2.2.4]{w1} for a proof.
\end{Proof}

\begin{Definition} \label{rdwaveletdefx}
For an integer $n$, let $X_n := X_{\iota_n}$ and let $\hat{X_n} := \hat{X_{\iota_n}}$ (see Definition \ref{rddefx}). Equip $X_n$ with the $n$th level bracket product $[ \; , \: ]_n$ and $n$th level module action $\circ_n$. We call $X_n$ the {\em $n$th level wavelet Hilbert $C^*(\BZ^d)$-module}, and $(X_n)_{n \in \BZ}$ a {\em wavelet chain of Hilbert $C^*(\BZ^d)$-modules}.

Equip $\hat{X_n}$ with the Fourier transformed $n$th level bracket product $\hl \; , \: \hr_n$ and $n$th level Fourier transformed module action $\widehat{\circ_n}$. We call $\hat{X_n}$ the {\em  $n$th level Fourier transformed wavelet Hilbert $C^*(\BZ^d)$-module}, and $(X_n)_{n \in \BZ}$ a {\em Fourier transformed wavelet chain of Hilbert $C^*(\BZ^d)$-modules}.
\end{Definition}

\begin{Theorem}
\label{rdwavelethmexists}
Let $\tilde{\dil} \in M^d(\BZ)$ be a dilation matrix and consider the standard multiresolution structure with index $m$ on $L^2(\BR^d)$ corresponding to $\tilde{\dil}$ with translation group $\Gamma = \BZ^d$ and dilation $\dil$. Let $n$ be an integer and let $X_n \subset L^2(\BR^d)$ and $\hat{X_n} \subset L^2(\widehat{\BR^d})$ be defined as in Definition \ref{rdwaveletdefx}. We have that
\begin{enumerate}
\item The space $C_c(\BR^d)$ is an inner product $C_c(\BZ^d)$-module with operations $[ \;, \: ]_n$ and $\circ_n$, and $\| \cdot \|_{X_n}$ is equal to the Hilbert module norm on $C_c(\BR^d)$.
\item The space $X_n$ is a full Hilbert $C^*(\BZ^d)$-module with the inner product $[ \; , \:]_n$. The module action of $C^*(\BZ^d)$ on $X_n$ has the property that for $f \in X_n$ and $a \in l^1(\BZ^d)$ (note that $l^1(\BZ^d) \subset C^*(\BZ^d)$), the module action of $a$ on $f$ is given by $f \circ_n a$. We furthermore have that if $(f_n)_{n \in \BN}$ is a convergent sequence in $X_n$, then it is a convergent sequence in $L^2(\BR^d)$. 
\item The operations $\hl \; , \: \hr_n$ and $\widehat{\circ_n}$ make $\hat{X_n}$ into a full Hilbert $C(\BT^d)$-module.
\item If we identify the isomorphic $C^*$-algebras $C^*(\BZ^d)$ and $C(\BT^d)$, then the Fourier transform on $\BR^d$ is a unitary isomorphism between $X_n$ and $\hat{X_n}$. 
\end{enumerate}
\end{Theorem}
\begin{Proof}{}
Part 1 is a direct result of Proposition \ref{thetainnerprodmod}. Parts 2,3 and 4 follow directly from Theorem \ref{rdhmexists}.
\end{Proof}

In Definition \ref{def:adjoint}, we introduced adjointable operators, the main morphisms between Hilbert modules. The following corollary demonstrates that  the translations and dilations are unitary adjointable operators. Recall from Definition \ref{def:isohm} that two Hilbert modules are isomorphic if there is a unitary adjointable operator from one to the other. 
\begin{Corollary}
\label{tdadjointable}
For $\gamma \in \BZ^d$ and integer $n$, the translation $\pi^n_\gamma$ (when restricted to $X_n$) is a unitary element of $\mathcal{L}(X_n)$. The dilation $\dil$ is a unitary element of $\mathcal{L}(X_n, X_{n+1})$, and $\hat{\dil}$ is a unitary element of $\mathcal{L}(\hat{X_n}, \hat{X_{n+1}})$. Hence the Hilbert $C^*(\BZ^d)$-modules $(X_n)_{n \in \BZ}$ are all isomorphic (as Hilbert $C^*(\BZ^d)$-modules) to each other and the Hilbert $C(\BT^d)$-modules $(\hat{X_n})_{n \in \BZ}$ are also isomorphic to each other.
\end{Corollary}
\begin{Proof}{}
The translation $\pi^n_\gamma$ maps $X_n$ to itself. It is the case that for all $f, g \in X_n$, 
$$
[ \pi^n_\gamma f, g]_n = [f, (\pi^n_\gamma)^{-1} g]_n
$$
for all $\gamma \in \Gamma$, and so $\pi^n_\gamma$ is an adjointable operator on $X_n$. The translation satisfies $\gamma^* = \gamma^{-1}$ and hence is a unitary operator on Hilbert modules. 

We now verify that the dilation $\dil$ maps $X_n$ onto $X_{n+1}$ and that $\hat{\dil}$ maps $\hat{X_n}$ onto $\hat{X_{n+1}}$. By Definition \ref{rdwaveletdefx}, $X_n$ is the completion of $C_c(\BR^d)$ with respect to the norm $\| \cdot \|_{X_n}$ and $\hat{X_n}$ is the Fourier transform of $X_n$. For $f \in C_c(\BR^d)$, we have
\bean
\|f \|_{X_{n+1}} &=& \sup_{\xi \in \BT^d} | \hl \hat{f}, \hat{f} \hr_{n+1} (\xi) |^{\half} \\
&=& \sup_{\xi \in \BT^d} | \hl \dil^{-1} \hat{f}, \dil^{-1} \hat{f} \hr_{n} (\xi) |^{\half} \\
&=& \| \dil^{-1} f \|_{X_n}.
\eean
We therefore have that $\dil$ maps $X_n$ onto $X_{n+1}$. It follows that $\hat{\dil}$ maps $\hat{X_n}$ onto $\hat{X_{n+1}}$ because $\hat{\dil} = \mathcal{F}^* \dil \mathcal{F}$. 

From equation (\ref{dilident}) it follows that for $f \in X_n$, $g \in X_{n+1}$,
$$
{[\dil f, g]_{n+1}} = [f, \dil^{-1}g]_n
$$
so $\dil^{-1} = \dil^*$ (where in this case $\dil^*$ is the adjoint of $\dil$ as an adjointable operator between Hilbert modules). Hence $\dil$ is a unitary element of $\mathcal{L}(X_n, X_{n+1})$.
From equation (\ref{dilhatident}) we have that for $p \in \hat{X_n}$, $q \in \hat{X_{n+1}}$,
$$
\hl \hat{\dil} p, q \hr_{n+1} = \hl p, \hat{\dil}^{-1}q \hr_n
$$
so $\hat{\dil}^{-1} = \hat{\dil}^*$. Hence $\hat{\dil}$ is a unitary element of $\mathcal{L}(\hat{X_n}, \hat{X_{n+1}})$.
\end{Proof}

The next proposition demonstrates that each $X_n$ shares the same linear space and is similar to \cite[Proposition 1.11]{pr2}. The difference between each $X_n$ is therefore in how the $C^*(\BZ^d)$-valued inner product and the module action are defined within the linear space. 

Recall from page \pageref{def:gmra} that we defined $\Delta \in \mathrm{Hom}(\BZ^d)$ to be $\Delta(\gamma) = \dil^{-1} \gamma \dil$ (as a unitary operator on $\mathcal{H}$), where the translation $\gamma \in \BZ^d$ is thought of as a unitary operator on $\mathcal{H}$. We shall be interested in the dual homomorphism $\hat{\Delta} \in \mathrm{Hom}(\BT^d)$. This satisfies 
$$ \label{Deltadual}
(\Delta(\gamma), \zeta) = (\gamma, \hat{\Delta}(\zeta)), \ \ \ \gamma \in \BZ^d, \zeta \in \BT^d.
$$
It follows that if we think of $\BT^d$ as being the quotient $\widehat{\BR^d} / \BZ^d$, then 
\be
\hat{\Delta}(\zeta) = \iota ( \tilde{\dil}^*(\zeta)), \ \ \ \zeta \in \BT^d.
\ee
In the above equation $\iota$ is the quotient map from $\BR^d$ onto $\BT^d$.

\begin{Proposition} \label{xnthesame}
For $p, q \in \hat{X_n}$, we have that
\be \label{eq:xnsamecalc}
\hl p, q \hr_{n-1}(\zeta) = \frac{1}{m} \sum_{\hat{\Delta}(\omega) = \zeta} \hl p, q \hr_n (\omega).
\ee
This implies that
\be
m^{-1/2} \| p \|_{X_n} \leq \| p \|_{X_{n-1}} \leq \| p \|_{X_n}.
\ee
We therefore have that the Hilbert $C^*(\BZ^d)$-modules $(X_n)_{n \in \BZ}$ all share the same linear space, and that the Hilbert $C(\BT^d)$-modules $(\hat{X_n})_{n \in \BZ}$ also all share the same linear space.
\end{Proposition}
\begin{Proof}{}
Recall that for $f \in C_c(\BR^d)$, $[f,f]_{n-1}(\gamma) = [f,f]_n(\Delta(\gamma))$ for $\gamma \in \BZ^d$. We calculate
\bean
\hl p,q \hr_{n-1}(\zeta)
&=& m^{n-1} \sum_{\hat{\iota_{n-1}}(\xi) = \zeta} \overline{p(\xi)}q(\xi)  \\
&=& m^{n-1} \sum_{\hat{\Delta}(\hat{\iota_n}(\xi)) = \zeta} \overline{p(\xi)}q(\xi)  \\
&=& m^{n-1} \sum_{\hat{\Delta}(\omega) = \zeta} \sum_{\hat{\iota_n}(\xi) = \omega} \overline{p(\xi)}q(\xi) \\
&=& \frac{1}{m} \sum_{\hat{\Delta}(\omega) = \zeta} \hl p, q \hr_n (\omega).
\eean
Now let
$$
F := \{ \zeta \in \BT^d : \hat{\Delta}(\zeta) = 0 \}.
$$
Then $F$ is a subgroup of $\BT^d$ with $m$ elements (because the index of $\Delta(\BZ^d)$ in $\BZ^d$ is equal to $m$). Because $\hat{\Delta}$ is a homomorphism, it follows that if $\omega, \omega' \in \BT^d$ satisfy $\hat{\Delta}(\omega) = \hat{\Delta}(\omega')$, then there exists $\beta \in F$ such that $\omega = \omega' + \beta$. We therefore have that
\be
\hl p, q \hr_{n-1} (\zeta) = \frac{1}{m} \sum_{\beta \in F} \hl p,q \hr_n ( \tilde{\dil}^{*-1}(\zeta) + \beta).
\ee
Therefore 
\bean
\|p\|_{X_{n-1}} &=& \sqrt{\sup_{\zeta \in \BT^d} \hl p , p \hr_{n-1} (\zeta)} \\
&=& \sqrt{\sup_{\zeta \in \BT^d} \frac{1}{m} \sum_{\beta \in F} \hl p,q \hr_n ( \tilde{\dil}^{*-1}(\zeta) + \beta)} \\
&\leq&  \sqrt{\frac{1}{m} \sum_{\beta \in F} \sup_{\zeta \in \BT^d} \hl p,q \hr_n ( \tilde{\dil}^{*-1}(\zeta) + \beta)} \\
&=& \sqrt{\frac{1}{m} \sum_{\beta \in F} \|p\|_{X_n}}  \\
&=& \| p \|_{X_n}.
\eean
We now use equation (\ref{ftbpuseful}) and the fact that $(\tilde{\dil}^*)^n \BZ^d \subset (\tilde{\dil}^*)^{n+1} \BZ^d$ to obtain that
\bean
\hl p, p \hr_{n-1} (\zeta)
&=& m^{n-1} \sum_{\beta \in (\tilde{\dil}^*)^{n-1} \BZ^d} \overline{ p \left( (\tilde{\dil}^*)^{n-1} \zeta + \beta \right) } p \left( (\tilde{\dil}^*)^{n-1} \zeta + \beta \right) \\
&\geq& m^{n-1} \sum_{\beta \in (\tilde{\dil}^*)^{n} \BZ^d} \overline{ p \left( (\tilde{\dil}^*)^{n-1} \zeta + \beta \right) } p \left( (\tilde{\dil}^*)^{n-1} \zeta + \beta \right).
\eean
This implies that
$$
\hl p, p \hr_{n-1} \left( \iota (\tilde{\dil}^* \zeta) \right) \geq m^{-1} \hl p , p \hr_n (  \zeta )
$$
where $\iota$ is the quotient map from $\widehat{\BR^d}$ to $\BZ^d$.  Therefore
\bean
\|p\|_{X_{n-1}} &=& \sqrt{\sup_{\zeta \in \BT^d} \hl p , p \hr_{n-1} (\zeta) } \\
&=& \sqrt{\sup_{\zeta \in \BT^d} \hl p , p \hr_{n-1} \left( \iota (\tilde{\dil}^* \zeta) \right) } \\
&\geq& \sqrt{ \sup_{\zeta \in \BT^d} m^{-1} \hl p , p \hr_n ( \zeta ) } \\
&=& m^{-1/2} \| p \|_{X_n}.
\eean
It therefore follows that the norms $( \| \cdot \|_{X_n} )_{n \in \BZ}$ are all equivalent. Because the linear space $X_n$ is the completion of $C_c(\BR^d)$ with respect to the norm $ \| \cdot \|_{X_n}$, it follows that the Hilbert $C^*(\BZ^d)$-modules $(X_n)_{n \in \BZ}$ all share the same linear space, and that the Hilbert $C(\BT^d)$-modules $(\hat{X_n})_{n \in \BZ}$ also all share the same linear space.
\end{Proof}

\section{Filters and Hilbert $L^\infty(\BT^d)$-modules} \label{s5}

In this section we shall investigate filters by using the bracket product. We shall use this to motivate the construction of a chain of Hilbert $L^\infty(\BT^d)$-modules $Y_n$, for integer $n$. Each Hilbert module $Y_n$ has its linear space containing the linear space $\widehat{X_n}$ as a subspace, and uses the Fourier transformed bracket product for its $L^\infty(\BT^d)$-valued inner product.

\begin{Definition} \label{def:filters}
Let $\psi_1, \ldots \psi_{m-1}$ be a multiwavelet corresponding to a single scaling function $\varphi$. We define
\bea \label{eq:filters}
h &=& [\varphi, \dil^{-1} \varphi]_0  \\
g^i &=& [\varphi, \dil^{-1} \psi^i]_0 
\eea
for $i = 1 ,\ldots , m-1$. We call $h$ and $g^i$ the {\em scaling filter} and the {\em wavelet filters}.
\end{Definition}

In much of the literature, scaling filters are known as {\em low pass filters} and wavelet filters are known as {\em high pass filters}. The term {\em filter} is sometimes also used to refer to the Fourier transforms of $h$ and $g^i$. By Lemma \ref{dilcalc} it follows that $h = [\dil^n \varphi, \dil^{n-1} \varphi]_n $ and $g^i = [\dil^n \varphi, \dil^{n-1} \psi^i ]_n$. Because the multiresolution analysis and wavelet spaces associated with $\phi$ satisfy the property that $V_{n-1}$ and $W_{n-1}$ is contained in $V_n$, we can write
\bea \label{eq:scalingeq}
\dil^{n-1} \varphi =  \dil^n \varphi \circ_n [ \dil^n \varphi, \dil^{n-1} \varphi ]_n = \dil^n \varphi \circ_n h \label{eq:hms} \\
\dil^{n-1} \psi^i =  \dil^n \varphi \circ_n [ \dil^n \varphi, \dil^{n-1} \psi^i ]_n =\dil^n \varphi \circ_n g^i. \label{eq:hmw}
\eea
We shall now show that $\|h \|_2 = \| g^i \|_2 = 1$. We know that $\| \dil^{-1}\varphi \|_2 =1 $, and so by (\ref{eq:hms}), $\| \sum_{\gamma \in \BZ^d} h(\gamma) \varphi(\cdot - \gamma )\|_2 = 1$. Now because the translations of $\varphi$ by $\BZ^d$ form an orthonormal set, it follows from Pythagoras' theorem that $\sum_{\gamma \in \BZ^d} |h(\gamma)|^2 = \| \sum_{\gamma \in \BZ^d} h(\gamma) \varphi(\cdot - \gamma )\|_2$. Therefore $\|h\|_2 = \sum_{\gamma \in \BZ^d} |h(\gamma)|^2 =1$. The proof that $ \| g^i \|_2 = 1$ is the same.

Using the approach to filters described here, one can investigate things such as the fast wavelet transform and the cascade algorithm (\cite{w1}, Chapter 3). In \cite{w1}, Theorem 3.4.10, it is shown that the cascade algorithm converges in the topology of $X_0$. We shall now examine some examples of filters.

The Haar wavelet is given by $\psi = \chi_{[0,1/2)} - \chi_{[1/2,1)}$ where $\chi_S$ is the characteristic function of a set $S$. It has the corresponding scaling function $\varphi = \chi_{[0,1)}$. It can be shown using routine calculations that
\bean
h &=& [\varphi, \dil^{-1} \varphi ]_0 = 2^{-1/2} (e_0 + e_{-1}), \\
g &=& [\varphi, \dil^{-1} \psi ]_0 = 2^{-1/2} (e_0 - e_{-1})
\eean
where $e_i \in C_c(\BZ)$ is the element satisfying $e_i(i)=1$ and $e_i(j) = 0$  when $j \neq i$.

The Shannon wavelet is given by $\hat{\psi}(\xi) = e^{i \xi /2} \chi_{[-1,-1/2 ) \cup (1/2,1]}(\xi)$ with scaling function $\hat{\varphi} = \chi_{[-1/2,1/2)}$. Using Lemma \ref{lma:rdfouhilb}, we calculate
\bean
\hat{h}(\zeta) = \hl \hat{\varphi}, \widehat{\dil^{-1} \varphi} \hr_0 (\zeta)
&=& \sqrt{2} \sum_{k \in \BZ} \overline{\chi_{[-1/2, 1/2)}(\zeta + k)} \chi_{[-1/4,1/4)}(\zeta + k) \\
&=& \sqrt{2} \chi_{[-1/4, 1/4)}(\zeta).
\eean
And
\bean
\hat{g}(\zeta) = \hl \hat{\varphi}, \widehat{\dil^{-1} \psi} \hr_0 (\zeta)
&=& \sqrt{2} \sum_{k \in \BZ} \overline{\chi_{[-1/2, 1/2)}(\zeta + k)} e^{i \zeta} \chi_{[-1/2, - 1/4) \cup (1/4, 1/2]}(\zeta + k) \\
&=& \sqrt{2} e^{i \zeta} \chi_{[-1/2, - 1/4) \cup (1/4, 1/2]}.
\eean
We have that $\hat{h}$ and $\hat{g}$ are contained in $L^2(\BT) \cap L^\infty(\BT)$ but not $C(\BT)$. It thus follows from Proposition \ref{xnthesame} that $\psi \notin X_0$. This motivates us to define ``larger'' Hilbert modules.

Recall that $L^{\infty}(\BT^d)$ is the continuous dual space of $L^1(\BT^d)^*$ in the sense that elements of $L^\infty(\BT^d)$ are bounded linear functionals of elements of $L^1(\BT^d)$. And recall that if $(f_n)_{n \in \BN}$ is a sequence in $L^\infty(\BT^d)$, then $(f_n)_{n \in \BN}$ converges to $f$ in the {\em weak* topology} if for all $a \in L^1(\BT^d)$, $\lim_{n \to \infty} f_n(a)= f(a)$.
Let $\theta$ be an embedding of $\BZ^d$ in $\BR^d$ and let the linear space $Y_\theta$ be the set of measurable functions $p$ on $\widehat{\BR^d}$ for which the norm
\be \label{Ythetanorm}
\| p \|_{Y_\theta} := \sqrt{ \mathrm{ess} \: \mathrm{sup}_{\zeta \in \BT^d} \frac{1}{\det (A_\theta)} \sum_{\beta \in (A_\theta^*)^{-1} \BZ^d} \overline{p (\zeta + \beta)} p (\zeta + \beta)  	}
\ee
is finite. It has been shown in \cite{cola1} that the above series does not necessarily converge in norm in $L^\infty(\BT^d)$, so we instead require weak* convergence. If we compare this norm to $\|\cdot \|_{X_\theta}$, we see that the difference is that $\|\cdot \|_{Y_\theta}$ involves taking the essential supremum instead of the supremum and we are working in the Fourier domain. We can write
$$
\| p \|_{Y_\theta} = \sqrt{ \mathrm{ess} \: \mathrm{sup}_{\zeta \in \BT^d} \hl p , p \hr_\theta (\zeta)}.
$$

\begin{Lemma}
\label{lma:weakstar}
Suppose that $p \in Y_\theta$. Then the sum in $\hl p , p \hr_\theta$ converges in the weak* topology to an element of $L^\infty(\BT^d)$.
\end{Lemma}
\begin{Proof}{}
To simplify our notation we shall assume without any loss of generality that $\theta$ is the natural embedding of $\BZ^d$ in $\BR^d$. We therefore want to show that the sum defined for almost every $\zeta \in \BT^d$ by
\be \label{Ythetaonenorm}
\sum_{\beta \in \BZ^d} \overline{p (\zeta + \beta)} p (\zeta + \beta)  
\ee
converges in the weak* topology to an element of $L^\infty(\BT^d)$. Suppose that $(S_n)_{n \in \BN}$ is a sequence of finite subsets of $\BZ^d$ for which $S_n \subset S_{n+1}$ and $\cup_{n \in \BN} S_n = \BZ^d$. Then 
$$
\sum_{\beta \in \BZ^d} \overline{p (\zeta + \beta)} p (\zeta + \beta)
= \lim_{n \to \infty} \sum_{\beta \in S_n} \overline{p (\zeta + \beta)} p (\zeta + \beta).
$$
Note that the summands in the above sum are positive so if it converges, it converges unconditionally. Because $p \in Y_\theta$, it is the case that the set of all $\zeta \in \BT^d$ for which $\sum_{\beta \in \BZ^d} \overline{p (\zeta + \beta)} p (\zeta + \beta)$ does not converge has measure zero. Define a function $L$ on $\BT^d$ by
$$
L(\zeta) = \left\{ \begin{array}{rl}
\sum_{\beta \in \BZ^d} \overline{p (\zeta + \beta)} p (\zeta + \beta) & \mbox{when it converges,} \\
0 & \mbox{otherwise.} \end{array} \right.
$$
Because
$$
\mathrm{ess} \: \mathrm{sup}_{\zeta \in \BT^d} \sum_{\beta \in \BZ^d} \overline{p (\zeta + \beta)} p (\zeta + \beta)  < \infty,
$$
it follows that $L \in L^\infty(\BT^d)$.

We shall now show that the sum in (\ref{Ythetaonenorm}) converges to $L$ in the weak* topology. We want to show that for all $a \in L^1(\BT^d)$,
$$
\lim_{n \to \infty} \int_{\BT^d} \sum_{\beta \in S_n} \overline{p (\zeta + \beta)} p (\zeta + \beta)  a(\zeta) d \zeta
= \int_{\BT^d}  L(\zeta) a(\zeta) d \zeta.
$$
Because all of the terms in the above sum except for $a(\zeta)$ are positive, it is sufficient to show that the above relation holds for $|a(\zeta)|$ (instead of $a(\zeta)$). Suppose that $E$ is a measurable subset of $\BT^d$, let $\phi_a(E) = \int_E |a(\zeta)| d \zeta$, then by \cite[Theorem 1.29]{ru1}, $\phi_a$ is a measure and for all measurable $g: \BT^d \to [0, \infty ]$,
$$
\int_{\BT^d} g(\zeta) d \phi_a(\zeta) = \int_{\BT^d} g(\zeta) a(\zeta) d \zeta .
$$
It now follows from Lebesgue's monotone convergence theorem (see \cite[Theorem 1.26]{ru1}) that 
$$
\lim_{n \to \infty} \int_{\BT^d} \sum_{\beta \in S_n} \overline{p (\zeta + \beta)} p (\zeta + \beta)  d \phi_a(\zeta)
= \int_{\BT^d}  L(\zeta) d \phi_a(\zeta).
$$
This verifies that the sum (\ref{Ythetaonenorm}) converges to $L$ in the weak* topology.
\end{Proof}

We note that it follows from the polarisation identity that if $p, q \in Y_\theta$, then $\hl p , q \hr_\theta \in L^\infty(\BT^d)$. For $p \in Y_\theta$,
$$
\| p \|_2^2 = \int_{\BT^d} \hl p,p \hr_{\theta} (\zeta) d \zeta \leq \mathrm{ess} \: \mathrm{sup}_{\zeta \in \BT^d} \hl p,p \hr_{\theta}(\zeta) = \|p \|_{Y_\theta}^2.
$$
and so $Y_\theta$ is contained in $L^2(\widehat{\BR^d})$.

\begin{Theorem}
\label{thm:linfty}
The space $Y_\theta$ is a Hilbert $L^\infty(\BT^d)$-module when equipped with the Fourier transformed bracket product $\hl \;, \: \hr_\theta$ with convergence in the weak* topology on $L^\infty(\BT^d)$, and associated module action $\widehat{\circ_\theta}$.
\end{Theorem}
\begin{Proof}{}
The space $Y_\theta$ satisfies Properties 1,2 and 3 of Definition \ref{def:hmodule} by Lemma \ref{rdftbracident}. For $p \in Y_\theta$, the sum
$$
\sum_{\beta \in (A_\theta^*)^{-1} \BZ^d} \overline{p(\zeta + \beta)} p(\zeta + \beta)
$$
is nonnegative for all $\zeta \in \BT^d$. Therefore $\hl p, p \hr_\theta$ is a positive element of $L^\infty(\BT^d)$, verifying property 4 of Definition \ref{def:hmodule}. If $\hl p,p \hr_\theta = 0$, then $\| p \|_{Y_\theta} = 0$ and hence $\| p \|_2 = 0$ since $\| p \|_2 \leq \|p \|_{Y_\theta}$. This implies that $p = 0$, verifying property 5 of Definition \ref{def:hmodule}.

We now show that the space $Y_\theta$ is complete. Suppose that $(p_j)$ is a Cauchy sequence in $Y_\theta$. Then for all $\varepsilon > 0$, there exists a natural number $J$ such that 
$$
j,k > J \Rightarrow \| p_j - p_k \|^2_{Y_\theta} < \varepsilon.
$$
This implies that for almost every $\zeta \in \BT^d$,
\bean
\hl p_j - p_k , p_j - p_k \hr_\theta (\zeta) &<& \varepsilon \\
\mbox{i.e.  } \sum_{\beta \in (A_\theta^*)^{-1} \BZ^d} \overline{ (p_j - p_k) (\zeta + \beta)} (p_j - p_k)  (\zeta + \beta)  &<& \varepsilon \det A_\theta \\
\mbox{so for all $\beta \in (A_\theta^*)^{-1} \BZ^d$,  } \ \ \ \overline{(p_j - p_k)  (\zeta + \beta)} (p_j - p_k)  (\zeta + \beta)  &<& \varepsilon \det A_\theta .
\eean
So for almost every $\xi \in \BR^d$, there exists a scalar $p(\xi)$ for which $\lim_{j \to \infty} p_j(\xi) = p(\xi)$. This defines a function $p$ almost everywhere on $\BR^d$.

We shall now show that $\lim_{k \to \infty} \|p - p_k \|_{Y_\theta} = 0$. Suppose that $(S_n)_{n \in \BN}$ is a sequence of finite subsets of $(A_\theta^*)^{-1} \BZ^d$ for which $S_n \subset S_{n+1}$ and $\cup_{n \in \BN} S_n = (A_\theta^*)^{-1} \BZ^d$. Now $p - p_k = \lim_{j \to \infty} p_j - p_k$, so
\bean
\|p - p_k \|_{Y_\theta}
&=& \sqrt{ \mathrm{ess} \: \mathrm{sup}_{\zeta \in \BT^d} \frac{1}{\det A_\theta} \sum_{\beta \in (A_\theta^*)^{-1} \BZ^d}
\lim_{j \to \infty} \overline{ (p_j - p_k) (\zeta + \beta)} (p_j - p_k)  (\zeta + \beta) } \\
&=& \sqrt{ \mathrm{ess} \: \mathrm{sup}_{\zeta \in \BT^d} \frac{1}{\det A_\theta} \lim_{n \to \infty} \sum_{\beta \in S_n}
\lim_{j \to \infty} \overline{ (p_j - p_k) (\zeta + \beta)} (p_j - p_k)  (\zeta + \beta) } \\
&=& \lim_{j \to \infty} \sqrt{ \mathrm{ess} \: \mathrm{sup}_{\zeta \in \BT^d} \frac{1}{\det A_\theta}  \lim_{n \to \infty} \sum_{\beta \in S_n} \overline{ (p_j - p_k) (\zeta + \beta)} (p_j - p_k)  (\zeta + \beta) } \\
&=& \lim_{j \to \infty} \sqrt{ \mathrm{ess} \: \mathrm{sup}_{\zeta \in \BT^d} \frac{1}{\det A_\theta} \sum_{\beta \in (A_\theta^*)^{-1} \BZ^d} \overline{ (p_j - p_k) (\zeta + \beta)} (p_j - p_k)  (\zeta + \beta) }
\eean
Because $(p_j)$ is Cauchy in $Y_\theta$, it follows that for all $\varepsilon > 0$, there exists a natural number $k$ such that $\|p - p_k \|_{Y_\theta} < \varepsilon$. Since $p = (p - p_k) + p_k$, it follows that $p \in Y_\theta$. Therefore $Y_\theta$ is complete. 
\end{Proof}

\begin{Definition} \label{def:linftyn}
For integer $n$, let $Y_n$ be the set of measurable functions $p$ on $\widehat{\BR^d}$ for which the norm
\be \label{eq:linftyn}
\| p \|_{Y_n} = \sqrt{ \mathrm{ess} \: \mathrm{sup}_{\zeta \in \BT^d} \hl p,p \hr_{n}(\zeta) }
\ee
is finite. Again we only require weak* convergence in the above series.  Equip $Y_n$ with the $n$th level Fourier transformed bracket product and the associated module action. By Theorem \ref{thm:linfty} it is a Hilbert $L^\infty(\BT^d)$-module. Note that by definition $\widehat{X_n} \subset Y_n$. We call $(Y_n)_{n \in \BZ}$ a {\em wavelet chain of Hilbert $L^\infty(\BT^d)$-modules}.
\end{Definition}

Because $p \in Y_\theta$ implies $\|p\|_2 \leq \|p\|_{Y_\theta}$, $p \in Y_n$ implies $\|p\|_2 \leq \|p\|_{Y_n}$.

The following proposition gives necessary and sufficient conditions for a set of elements of $L^2(\BR^d)$ to be a multiwavelet. It also demonstrates that any wavelet will have it's Fourier transform contained in $Y_0$. 

\begin{Proposition}
\label{neccsuffnicey}
Suppose that $(Y_n)_{n \in \BZ}$ is a wavelet chain of Hilbert $L^\infty(\BT^d)$-modules. Suppose that $\psi^1, \ldots, \psi^M$ are elements of $L^2(\BR^d)$. Then 
$\{ \psi^1, \ldots, \psi^M \}$ is an orthonormal multiwavelet if and only if
\begin{enumerate}
\item For $i,j = 1, \ldots, M$, and all integers $m,n$,
\be \label{eq:onwc1y}
{\hl \hat{\dil^n} \hat{\psi^i}, \hat{\dil^m} \hat{\psi^j} \hr_n} = \delta_{i,j} \delta_{m,n} \mathbf{1}
\ee

where $\delta_{i,j}$ is the Kronecker delta and $\mathbf{1}$ is the function on $\BT^d$ which takes the value $1$ everywhere.

\item For all $f \in L^2(\BR^d)$,
\be \label{eq:onwc2y}
\hat{f} = \sum_{i = 1}^M \sum_{n \in \BZ} \hat{\dil^n} \hat{\psi^i} \widehat{\circ}_n \hl \hat{\dil^n} \hat{\psi^i}, \hat{f} \hr_n
\ee

and the above sum converges in $L^2(\BR^d)$.
\end{enumerate}
Furthermore, if $\{ \psi^1, \ldots, \psi^M  \}$ is a multiwavelet, then for all $i$, $\hat{\psi^i} \in Y_0$. If $\{ \varphi^1, \dots, \varphi^r \}$ is a set of scaling functions, then for all $i$, $\hat{\varphi^i} \in Y_0$.
\end{Proposition}
\begin{Proof}{}
Suppose that $\{ \psi^1, \ldots, \psi^M \}$ is an orthonormal multiwavelet. Then from orthogonality it holds that for each $i=1, \ldots, M$, $\hl \hat{\psi^i}, 
\hat{\psi^i} \hr_0 = \mathbf{1}$ and so $\hl \widehat{\dil^n \psi^i}, \widehat{\dil^n \psi^i} \hr_n = \mathbf{1}$ for all integers $n$. This implies that $\| \widehat{\dil^n \psi^i} \|_{Y_n} = 1$.
We therefore have that $\widehat{\dil^n \psi^i} \in Y_n$ and in particular $\hat{\psi^i} \in Y_0$. We have that $\hat{\varphi^i} \in Y_0$, because $\hl \hat{\varphi^i}, \hat{\varphi^i} \hr_0 = \mathbf{1}$.
Equation (\ref{eq:onwc1y}) also follows from orthogonality. Equation (\ref{eq:onwc2y}) is a direct consequence of the 
fact that the wavelets form an orthonormal basis for $L^2(\BR^d)$, with the sum converging because the $L^2$-norm of $f$ is finite.

Now suppose that equations (\ref{eq:onwc1y}) and (\ref{eq:onwc2y}) hold. Then equation (\ref{eq:onwc1y}) implies that 
$$
\{ \dil^n(\gamma(\psi^i)) \}_{\gamma \in \BZ^d, n \in \BZ, i = 1, \ldots, M}
$$
is an orthonormal set and equation (\ref{eq:onwc2y}) implies that
$$\{ \dil^n(\gamma(\psi^i)) \}_{\gamma \in \BZ^d, n \in \BZ, i = 1, \ldots, M}
$$
is an orthonormal basis for $L^2(\BR^d)$. It therefore follows that $\{ \psi^1, \ldots, \psi^M \}$ is an orthonormal multiwavelet.
\end{Proof}

\section*{Acknowledgments}

I would like to thank Peter Dodds, Bill Moran, Jaroslav Kautsky, Adam Rennie, Fyodor Sukochev and the referees for useful suggestions and comments. I would also like to thank Professor Marc Rieffel for sending me a copy of the notes for his talk \cite{r3}.

This research was partially financially supported by the Australian government through an Australian Postgraduate Award scholarship, and by the Cooperative Research Centre for Sensor, Signal and Information Processing (CSSIP) with a CSSIP supplementary scholarship.


{\footnotesize
\centerline{\rule{9pc}{.01in}}
\medskip

\centerline{School of Informatics and Engineering}
\centerline{Flinders University of South Australia}
\centerline{GPO Box 2100, Adelaide, SA, Australia 5001}
\centerline{email: pwood@infoeng.flinders.edu.au, pjwood@myplace.net.au} 
}
\end{document}